\documentclass[twocolumn, switch]{article} 

\usepackage{preprint}
\usepackage{subcaption}

\usepackage{amsmath, amsthm, amssymb, amsfonts}

\usepackage[numbers,square]{natbib}
\usepackage[utf8]{inputenc}	
\usepackage[T1]{fontenc}	
\usepackage{xcolor}		
\usepackage[colorlinks = true,
            linkcolor = purple,
            urlcolor  = blue,
            citecolor = cyan,
            anchorcolor = black]{hyperref}	
\usepackage{booktabs} 		
\usepackage{nicefrac}		
\usepackage{microtype}		
\usepackage{lineno}		
\usepackage{float}			

\usepackage{lipsum}		

\usepackage{tabularray}

\usepackage{nomencl}
\makenomenclature

\usepackage{newfloat}
\DeclareFloatingEnvironment[name={Supplementary Figure}]{suppfigure}
\usepackage{sidecap}
\sidecaptionvpos{figure}{c}

\usepackage{titlesec}
\titlespacing\section{0pt}{12pt plus 3pt minus 3pt}{1pt plus 1pt minus 1pt}
\titlespacing\subsection{0pt}{10pt plus 3pt minus 3pt}{1pt plus 1pt minus 1pt}
\titlespacing\subsubsection{0pt}{8pt plus 3pt minus 3pt}{1pt plus 1pt minus 1pt}

\usepackage{tikz}

\definecolor{lime}{HTML}{A6CE39}
\DeclareRobustCommand{\orcidicon}{
	\begin{tikzpicture}
	\draw[lime, fill=lime] (0,0) 
	circle [radius=0.16] 
	node[white] {{\fontfamily{qag}\selectfont \tiny ID}};
	\draw[white, fill=white] (-0.0625,0.095) 
	circle [radius=0.007];
	\end{tikzpicture}
	\hspace{-2mm}
}
\foreach \x in {A, ..., Z}{\expandafter\xdef\csname orcid\x\endcsname{\noexpand\href{https://orcid.org/\csname orcidauthor\x\endcsname}
			{\noexpand\orcidicon}}
}

\newcommand{\G}{\mathbf{G}}
\newcommand{\N}{\mathbf{N}}
\newcommand{\A}{\mathbf{A}}
\newcommand{\R}{\mathbb{R}}
\newcommand{\Z}{\mathbb{Z}}
\newcommand{\E}{\mathbb{E}}
\newcommand{\PR}{\mathbb{P}}
\newcommand{\x}{\mathbf{x}}
\newcommand{\bu}{\mathbf{u}}
\newcommand{\z}{\mathbf{z}}
\newcommand{\ba}{\mathbf{a}}
\newcommand{\bb}{\mathbf{b}}
\newcommand{\bs}{\mathbf{S}}
\newcommand{\T}{\mathbf{T}}
\newcommand{\J}{\mathbf{J}}
\newtheorem{definition}{Definition}[section]

\usepackage{algorithm,algorithmic}
\newcommand\algorithmicprocedure{\textbf{procedure}}
\newcommand{\algorithmicendprocedure}{\algorithmicend\ \algorithmicprocedure}
\makeatletter
\renewcommand{\ALG@name}{Pseudo-Code}
\newcommand\PROCEDURE[3][default]{%
  \ALC@it
  \algorithmicprocedure\ \textsc{#2}(#3)%
  \begin{ALC@prc}%
}
\newcommand\ENDPROCEDURE{%
  \end{ALC@prc}%
  \ifthenelse{\boolean{ALC@noend}}{}{%
    \ALC@it\algorithmicendprocedure
  }%
}
\newenvironment{ALC@prc}{\begin{ALC@g}}{\end{ALC@g}}
\makeatother

\usepackage{enumitem}

\title{An Integrated Cyber-Physical Risk Assessment Framework for Worst-Case Attacks in Industrial Control Systems\thanks{This work was funded by National Science Foundation Award 2119654.}}

\usepackage{eso-pic,picture}

\AddToShipoutPicture{%
\setlength{\unitlength}{1pt}%
\protect{\put(.5\paperwidth,.5\paperheight){%
\makebox(0,-745){\rotatebox{0}{\textcolor[gray]{0.60}%
   {\color{gray}{*correspondence: \texttt{naftabi@clemson.edu}}}}}
\makebox(0,755){\rotatebox{0}{\footnotesize\textcolor[gray]{0.60}%
   {\color{gray}{This work has been submitted to the IEEE for possible publication.}}}}

\makebox(0,735){\rotatebox{0}{\footnotesize\textcolor[gray]{0.60}%
   {\color{gray}{Copyright may be transferred without notice, after which this version may no longer be accessible.}}}}
   
   }}\ClearShipoutPictureBG}

\usepackage{abstract}
\usepackage[affil-it]{authblk}

\author[1\thanks{\tt{naftabi@clemson.edu}}]{Navid Aftabi}
\author[1\thanks{\tt{dli4@clemson.edu}}]{Dan Li, Ph.D.}
\author[1\thanks{\tt{tcshark@clemson.edu}}]{Thomas Sharkey, Ph.D.}

\affil[1]{Department of Industrial Engineering, Clemson University, Clemson, SC 29634, USA}

\begin{document}

\twocolumn[ 
\begin{@twocolumnfalse} 
  
\maketitle

\begin{abstract}
Industrial Control Systems (ICSs) are widely used in critical infrastructures that face various cyberattacks causing physical damage. With the increasing integration of the ICSs and information technology (IT), ensuring the security of ICSs is of paramount importance. In an ICS, cyberattacks exploit vulnerabilities to compromise sensors and controllers, aiming to cause physical damage. Maliciously accessing different components poses varying risks, highlighting the importance of identifying high-risk cyberattacks. This aids in designing effective detection schemes and mitigation strategies. This paper proposes an optimization-based cyber-risk assessment framework that integrates cyber and physical systems of ICSs. The framework models cyberattacks with varying expertise and knowledge by 1) maximizing physical impact in terms of time to failure of the physical system, 2) quickly accessing the sensors and controllers in the cyber system while exploiting limited vulnerabilities, 3) avoiding detection in the physical system, and 4) complying with the cyber and physical restrictions. These objectives enable us to model the interactions between the cyber and physical systems jointly and study the critical cyberattacks that cause the highest impact on the physical system under certain resource constraints. Our framework serves as a tool to understand the vulnerabilities of an ICS with a holistic consideration of cyber and physical systems and their interactions and assess the risk of existing detection schemes by generating the worst-case attack strategies. We illustrate and verify the effectiveness of our proposed method in a numerical and a case study. The results show that a worst-case strategic attacker causes almost 19\% further acceleration in the time to failure of the physical system while remaining undetected compared to a random attacker. 
\end{abstract}
\keywords{Cyber-Physical Systems \and Cybersecurity \and Industrial Control Systems \and Worst-Case Attack \and Optimization \and Reliability Engineering \and Risk Assessment} 
\vspace{0.35cm}

\end{@twocolumnfalse} 
] 
\saythanks


\nomenclature{ICS}{Industrial Control Systems}
\nomenclature{IT}{Information Technology}
\nomenclature{MTTF}{Mean Time to Failure}
\nomenclature{BM}{Brownian Motion}
\nomenclature{MILP}{Mixed-Integer Linear Program}
\nomenclature{HIL}{Hardware-in-the-Loop}
\nomenclature{BWPP}{Boiling Water Power Plant}
\nomenclature{WV}{Water Valve}
\nomenclature{SV}{Steam Valve}
\nomenclature{FV}{Fuel Valve}
\nomenclature{AH}{Administrator Host}
\nomenclature{HM}{Human-Machine Interface}
\nomenclature{DS}{Data Server}
\nomenclature{PCs}{Personal Computers}
\nomenclature{DCS}{Distributed Control System}
\nomenclature{s.t.}{Subject To}
\nomenclature{o.w.}{Otherwise}
\nomenclature{TTF}{Time to Failure}

\printnomenclature[0.6in]


\section{Introduction}
\begin{figure}[!t]
    \centering
    \includegraphics[width=.9\linewidth,page=1]{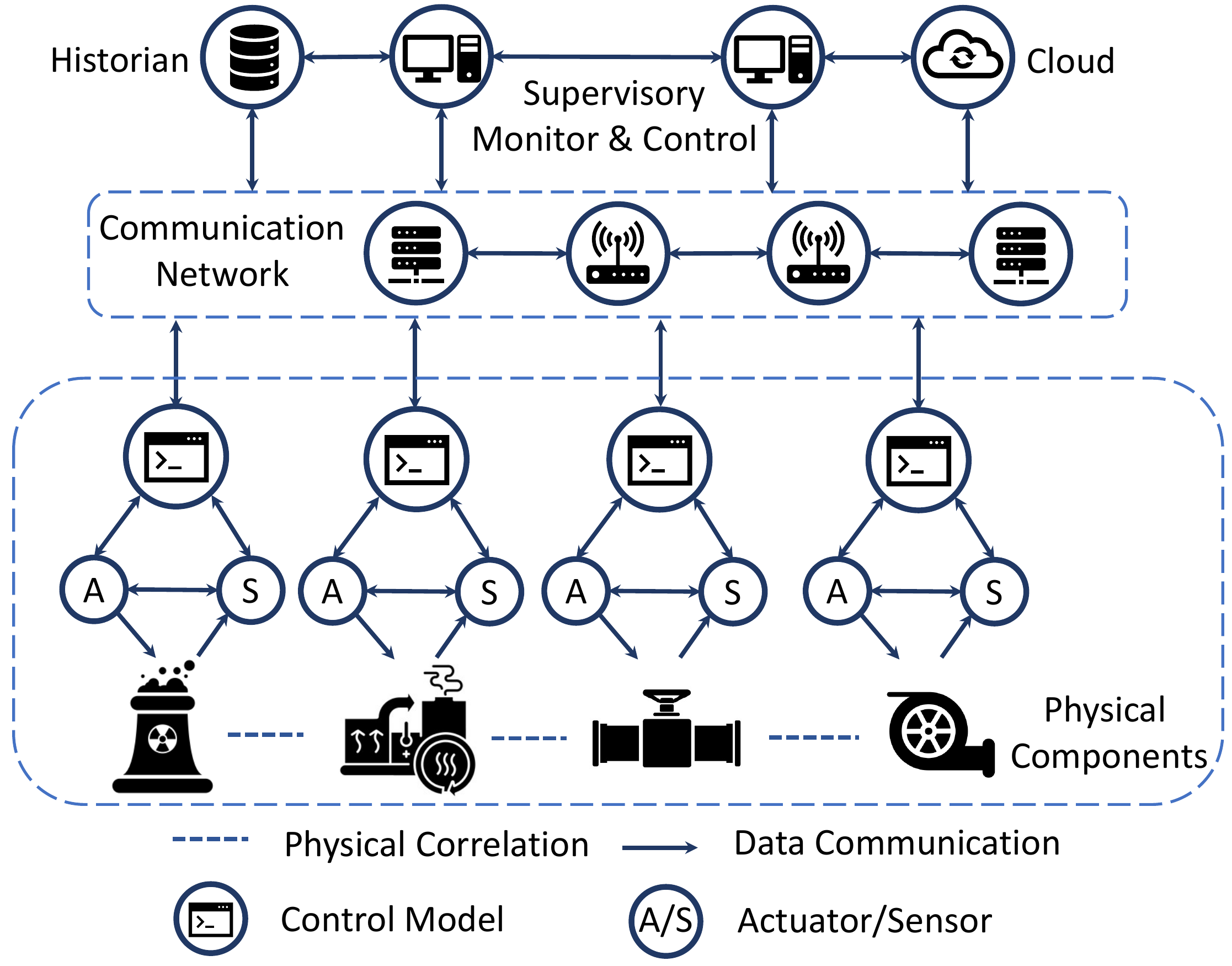}
    \caption{Demonstration of a generic ICS network structure.}
    \label{fig:cps}
\end{figure}

In recent decades, cybersecurity for ICSs has become increasingly critical. A typical multi-component ICS (as shown in Figure~\ref{fig:cps}) consists of multiple cyber components, such as a communication network, supervisory monitor and control, cloud, and historian, where heterogeneous components interact physically and are linked through data communication between supervisory and local control networks while monitored by a distributed control system (DCS). These systems are crucial for automating and managing critical infrastructure like power and water systems. Originally designed for isolated networks, they have evolved with the Internet of Things, edge and cloud computing, and 5G technologies, making them smarter but more exposed to cyberattacks~\cite{cheminod2012review,zhang2022advancements}. Cyberattacks target this system by exploiting vulnerabilities to access critical sensors and controllers. Significant physical consequences of such attacks, like the 2009 Stuxnet and the 2015 Ukrainian power grid attacks, demonstrate the vital need for robust, resilient ICSs~\cite{ning2021design,zhang2022advancements}. Addressing this requires a comprehensive approach to cyber-risk management, encompassing risk assessment, mitigation strategies, and continuous security monitoring~\cite{cheminod2012review,khalil2023threat}.

Risk is generally the likelihood of an event adversely affecting an organization's objectives~\cite{cheminod2012review}. In information security, cyber-risk is the potential physical damage from threats exploiting vulnerabilities in assets~\cite{cheminod2012review}. A key challenge in ICSs is scarce attack data, leading to the use of threat modeling for cyber-risk assessment, a systematic process for pinpointing potential threats to a system~\cite{khalil2023threat}. A comprehensive ICS cyber-risk assessment needs a generalized threat modeling approach that holistically captures the structure of ICSs, particularly as cyberattacks lead to physical impact. The complexity and heterogeneity of ICSs also demand a multi-paradigm framework for effective cyber-risk assessment~\cite{barivsic2022multi}. Moreover, since cyberattack success is contingent upon both the system vulnerabilities and the attackers' ability to utilize them, the attackers with certain levels of knowledge and expertise on both aspects of ICSs strategically allocate their resources in targeting an ICS~\cite{cheminod2012review}. Consequently, studying the worst-case cyberattacks and the risk associated with them is crucial not only for strengthening the security of ICSs but also for enhancing their resilience through the potential development of effective detection mechanisms. Given the challenges outlined, several critical research questions emerge.

\textit{First, how can we develop a generalized threat modeling approach that comprehensively captures the intricate structure of ICSs? }
Existing threat modeling approaches often fail to consider the structure of the ICSs holistically (see the review of~\cite{khalil2023threat}). These methods tend to focus either solely on the cyber aspect~\cite{lallie2020review} or exclusively on the physical components~\cite{guo2016optimal,guo2018worst,jeon2019stealthy,zhang2019optimal,shang2021worst,zhang2021stealthy,biehler2023sage,anand2023risk}. However, the interplay between cyber and physical systems is critical in ICS operations, particularly since cyberattacks are aimed at causing physical damage. To address these limitations, an integrated cyber-physical model of the ICS is necessary to develop a generalized threat model of a cyberattack.

\textit{Second, how can cyberattacks be strategically executed to maximize their impact on the physical system? }
In implementing a cyberattack on an ICS, attackers often adopt a strategic approach, particularly when they possess sufficient knowledge of the system's operations. Current cyber-risk assessment methods, however, have a limited scope. Many either focus just on the cyber aspects of ICSs~\cite{enayaty2020survey, nguyen2017computational} without fully considering the physical consequences or concentrate on the physical systems while overlooking the cyber network's connectivity and vulnerabilities~\cite{teixeira2015secure, anand2023risk}. Although a significant part of existing literature includes cross-layer and cyber-to-physical risk assessments~\cite{huang2018assessing, huang2019game,cao2019operational,qin2020association, yang2021cross, deng2022quantitative}, these often assume a cyberattack and assess its impact by propagating the success likelihood from the cyber system to the physical system. They often do not capture the full complexity of attack strategies. A comprehensive threat model should include a physics-based objective as the cyber-to-physical impact that guides attackers in identifying and exploiting the most vulnerable points in the ICS. Such an objective leads to attack strategies that not only breach cybersecurity but also cause maximum physical disruption. To bridge this gap, it is required to develop a physics-based objective for a general attacker to prioritize their actions in targeting those system components whose compromise would lead to the highest physical damage.

\textit{Lastly, how can the worst-case attack with the highest physical impact be identified?} In the context of ICSs, such `worst-case attacks' are defined as cyberattacks that aim to cause maximum disruption or damage to the physical system by exploiting the cyber vulnerabilities and manipulating operational data under the limited knowledge and resources of the entire ICS. Modeling and analyzing extreme events play an important role in risk assessment~\cite{mcneil1999extreme}. Identifying the `worst-case attacks' provides valuable insights for the design of mitigation and detection schemes. To identify these attacks, an optimization framework needs to be developed to consider the generalized attack model, the physics-based objective, and varying levels of attacker knowledge and resources. However, the probabilistic approaches in the literature~\cite{huang2018assessing, huang2019game,cao2019operational,qin2020association, yang2021cross, deng2022quantitative} often rely on specific designs of attacks to assess their risks, lacking a comprehensive and systematic consideration of all possible attack strategies and solving for the optimal solution to characterize the `worst-case attack' in a generalizable fashion.

This paper responds to these research questions by introducing a novel optimization framework. This framework is designed to identify the strategic worst-case cyberattack in an ICS. The aim of the attacker is to gain control over sensors or controllers with limited effort, thereby manipulating the physical system into an undesirable state without triggering existing detection mechanisms. We focus on cyberattacks that compromise the availability and integrity of operational data and control, as these can inflict physical damage by tampering with both cyber and physical systems of ICSs. To verify the effectiveness, flexibility, and applicability of our proposed framework, we solve the framework for instances of synthetic ICSs containing both cyber and physical systems and conduct a case study on a boiling water reactor power plant (BWPP) ICS (introduced in~\cite{huang2018assessing}). The key contributions of this paper are as follows:

\begin{enumerate}
    \item We introduce a generalized data integrity attack model under the joint cyber-physical system framework, integrating an attack graph representing the cyber layer and a linear system dynamics model representing the physical layer.
    \item We leverage the physical system dynamics of the ICS and its degradation to quantify the cyber-to-physical impact in terms of the time to failure (TTF). More specifically, we propose a physics-based objective function for the attacker based on the ICS's mean time to failure (MTTF) intended to be minimized for imposing physical impact.
    \item Using the generalized attack model and the physics-based objective function, we formulate an optimization problem to identify worst-case attacks with the highest physical impact under different levels of the attacker's knowledge and resources.
\end{enumerate}
 
The rest of this paper is structured as follows: Section~\ref{sec:literature} reviews related literature and outlines research gaps and contributions. Section~\ref{sec:background} establishes the necessary background for our framework. Section~\ref{sec:method} details the proposed framework. Section~\ref{sec:experiment} presents experimental results, a real-world case study, and a comparative analysis. The paper concludes with a discussion in Section~\ref{sec:conclusion}.

\section{Related Work}
\label{sec:literature}
In this section, we briefly review the related literature on the cyber-risk assessment of ICSs in both the cyber and physical systems and discuss the relevant work supporting formalizing the cyber-aware physical impact of a cyberattack on the ICSs.
\subsection{Traditional Cyber-Risk Assessment}
From an optimization modeling perspective, several methods have been used for threat modeling and formalizing the ICSs as a framework and use them for threat identification, vulnerability analysis, cyber-risk assessment, and risk mitigation \cite{sharkey2021search,enayaty2020survey}. Stochastic programming \cite{zheng2019budgeted,zheng2019robust,schmidt2021risk,zheng2019interdiction}, multi-objective optimization \cite{khouzani2019scalable}, and network interdiction and multi-level optimization \cite{zheng2019interdiction,baycik2018interdicting} are some of the methods have been applied in the literature. The main focus of these methods is to increase the effort needed for cyberattacks to be completed. In this line of research, the focus is on understanding the types of attacks that can potentially be conducted within the cyberinfrastructure of a physical system. However, they do not model the impact of how a successful attack will then change the operations of the physical system.
\subsection{Physics-based Risk Assessment}
Recent literature presents various threat modeling and cyber-risk assessment methods. Notably,  these include the identification of the worst-case linear attack strategy with $\chi^{2}$~\cite{guo2016optimal} or Kullback-Leibler divergence~\cite{shang2021worst,guo2018worst} stealthiness constraints while maximizing the remote estimation error covariance. Other approaches involve deception attacks with Kullback-Leibler divergence stealthiness constraints while maximizing a quadratic utility function~\cite{zhang2019optimal}, and pole-dynamics attacks with measurement deviation stealthiness constraints while maximizing the norm of system state~\cite{jeon2019stealthy}. Furthermore, studies have explored stealthy integrity attacks modeled as a closed-loop dynamical system with model-based stealthiness constraints~\cite{zhang2021stealthy}. Recent research shifts towards identifying generic worst-case attack strategies by incorporating a generic stealthiness constraint while minimizing the norm of deviation of a damage function based on system output~\cite{biehler2023sage} and a model-based stealthiness constraint while maximizing output-to-output gain~\cite{anand2023risk}. These approaches primarily focus on a physical system rather than the entire ICS and do not consider the strategic cyber intrusion process. Furthermore, they typically assume a specific type of attack and presume that attackers have access to all or some attack channels, which is an unrealistic assumption.
\subsection{Cross-Layer Risk Assessment}
A noticeable limitation of the physics-based and traditional cyber-risk assessment approaches is their tendency to independently model each aspect of an ICS. Recently, efforts have been made for cyber-risk assessment and establishing risk mitigation strategies while considering both cyber and physical systems. These methods include cross-layer and cyber-to-physical risk assessment methods. These multi-paradigm frameworks usually propagate cyberattack success probabilities from the cyber system to the physical system. The Bayesian network with stochastic hybrid systems \cite{huang2018assessing}, Markov chain with network theory \cite{cao2019operational}, Bayesian analysis with game theory \cite{huang2019game}, data-driven models with Bayesian network \cite{qin2020association}, fuzzy Bayesian network with the linear state-space model \cite{yang2021cross}, and  Bayesian analysis with epidemic model\cite{deng2022quantitative} have been integrated to capture the structural characteristics of both of the systems. These studies mainly focus on the disruptions of physical systems by the propagated success probabilities instead of considering the entire complexity of ICS and the strategic cyber attack.
\subsection{Degradation and Remaining Useful Life}
An essential step in the cyber-risk assessment of the ICSs is to evaluate the physical outcome of the malicious actions of a cyberattack on both the cyber and physical systems. The most traditional risk assessment approaches evaluate the potential loss of a cyberattack in the monetary value. Even though this approach is sensible, it requires assumptions on the values of assets of ICSs. An appropriate alternative for evaluating the physical outcome is to assess how cyberattacks accelerate system failure. Studying the degradation process of the physical system provides a tool to introduce a new security metric to quantify the risk associated with cyberattacks on the ICSs. The literature on the TTF, degradation, and the remaining useful life of the physical systems is extensive. For instance, \cite{li2020multi,li2021degradation,peng2018switching} modeled the degradation signal of a physical process as a drifted Brownian motion (BM) that provides insights about the physical system lifetime.
\begin{figure*}[!t]
    \centering
    \includegraphics[width=.9\linewidth,page=2]{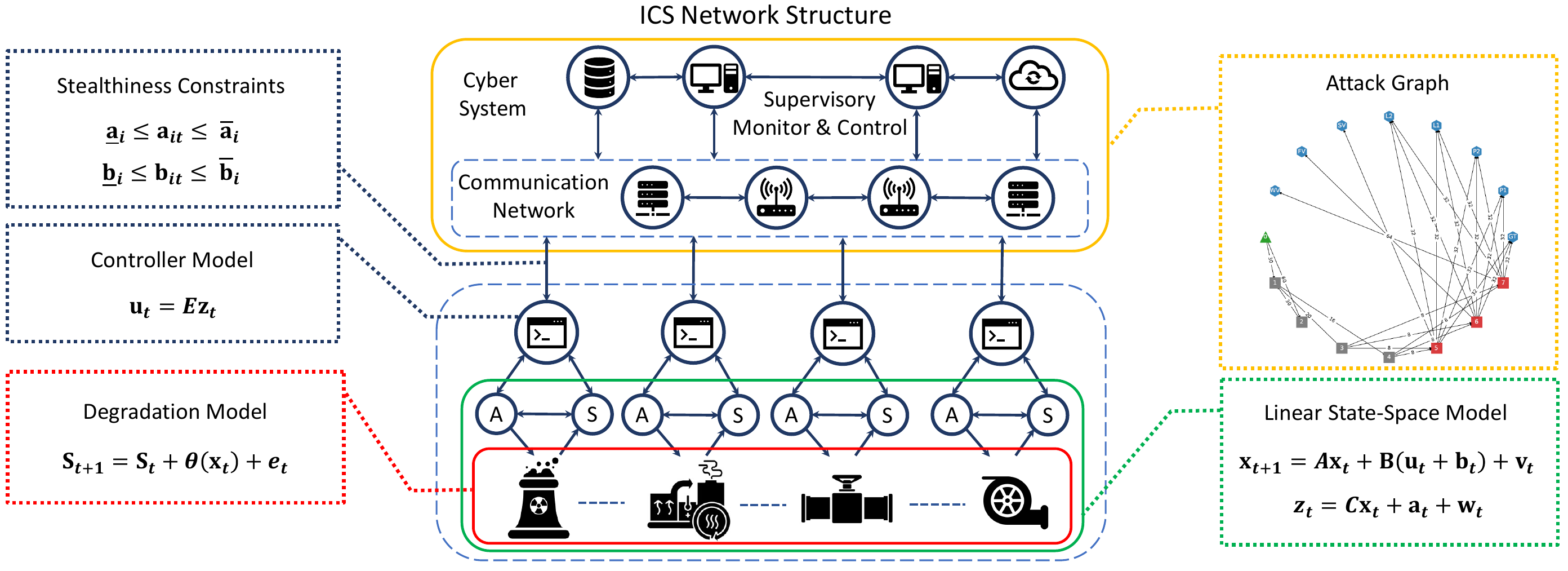}
    \caption{Proposed risk assessment framework.}
    \label{fig:fw}
\end{figure*}

\section{System Modeling}
\label{sec:background}
As discussed earlier, a general ICS consists of two interconnected networks, cyber and physical systems. This section describes the models and methods of characterizing each network and the cyber-to-physical outcomes.
\subsection{Cyber System Model}
\label{sec:cyber}
Among the attack modeling techniques, in the cyber system of an ICS, an attack graph or tree is one of the popular methods that mathematically and visually represents the sequence of events that lead to a successful cyberattack \cite{lallie2020review}. On this graph, nodes represent the attack states or security pre-/post-conditions, and arcs correspond to the transition of states fulfilled by an attack action. An attack state becomes active if its required attack state(s) are activated before. Consequently, a successful attack is a \textit{feasible path} from initial vulnerable states to a target state on the attack graph while satisfying the \textit{precondition dependencies}. A \textit{critical path} on this graph is a feasible path that requires the least effort. The effort or difficulty of vulnerability exploitation is usually represented by assigning a quantity on each arc. Depending on the range of values for this quantity, the cyber system can be represented as a probabilistic attack graph~\cite{khouzani2019scalable} or exploitation-time attack graph~\cite{zheng2019interdiction}.

Mathematically speaking, we define the attack graph as a directed network $\G=(\N,\A)$ where $\N$ and $\A$ are sets of its nodes and arcs representing attack states and exploitations, respectively. To reach an attack state $j\in\N$, the dependent preconditions $i\in\N$ must be satisfied such that $(i,j)\in\A$. Without loss of generality, we consider target nodes to be attack states where only a particular sensor or controller is being compromised. We define $\N_{s}\subset\N$ and $\N_{c}\subset\N$ as the set of states where sensors and controllers are compromised, respectively. Therefore, $|\N_{s}|$ and $|\N_{c}|$ equal the number of sensors and controllers in the physical system, respectively.
\subsection{Physical System Model}
\label{sec:physical}
The physical system being controlled is modeled as a multivariate discrete-time linear time-invariant system of equations. This formulation of the physical process allows us to: ($i$) quantify the cyber-to-physical impact, and ($ii$) it is easily generalizable to formulate many physical systems. We consider a multi-component physical system whose dynamics are modeled as
\begin{subequations}\label{eq:nsys}
    \begin{align}
    \x_{t+1} &= A\x_{t} + B\bu_{t}+\mathbf{v}_{t} \\
    \z_{t} &= D\x_{t} + \mathbf{w}_{t}
    \end{align}
\end{subequations}
In above equations, $\x_{t}\in\R^{|M|}$ represents the system states, and $M$ denotes the set of these states; $\bu_{t}\in\R^{|\N_{c}|}$ and $\z_{t}\in\R^{|\N_{s}|}$ are the control actions and sensor measurements respectively; $\mathbf{v}_{t}\in\R^{|M|}$ and $\mathbf{w}_{t}\in\R^{|\N_{s}|}$ are the process and measurement noises respectively; and, $A,B,$ and $C$ are matrices of appropriate dimensions.
\subsection{Degradation Model}
\label{sec:deg}
The physical system state variable $\x_{t}$ is used to model physical component degradation. The degradation signal $\bs(t)$ is typically modeled as a drifted BM, $\bs(t)=\mu t + \sigma_{s} \mathbf{B}(t)$ where $\mathbf{B}(t)$ is an standard BM such that $\mathbf{B}(t)\sim\mathcal{N}(0,t)$ that implies $\bs(t)\sim\mathcal{N}(\mu t, \sigma_{s}^{2}t)$. The degradation rate of the physical system is a function of the system state $\x_{t}$, representing the working condition of the system in practice. Given a sequence $\{\x_{t}\}$ and observing the value of the logged signal in the discrete time intervals, the discretized model of the degradation signal $\bs_{t}$ of the physical system is defined as
\begin{equation}\label{eq:degmodel}
    \bs_{t}=\bs_{t-1}+\theta(\x_{t})+ e_{t}
\end{equation}
where $e_{t}\sim\mathcal{N}(0,\sigma_{s}^{2})$ are independent and identically distributed Gaussian white noise for each $t$, and, the degradation rate $\theta(\x_{t})$ is a function of $\x_{t}$. To determine the distribution of $\bs_{t}$, we define the random variable $\J_{t}$ to be the jump in the degradation signal from time $t-1$ to $t$ for $t\geq1$, i.e., $\J_{t}=\bs_{t}-\bs_{t-1}$. Since $\{e_{t}\}$ is a sequence of independent random variables, $\{\J_{t}\}$ is also a sequence of independent random variables, such that $\J_{t}\sim\mathcal{N}(\theta(\mathbf{x}_{t}),\sigma^{2}_{s})$. Furthermore, the degradation signal can be written as the sum of independent jumps up to time $t$, i.e., $\bs_{t}=\sum_{\tau=1}^{t} \J_{\tau}$, that implies 
\begin{equation}\label{eq:dist}
\bs_{t}\sim\mathcal{N}\left(\sum\nolimits_{\tau=1}^{t} \theta(\mathbf{x}_{\tau}), \sigma^{2}_{s}t \right)
\end{equation}

\section{ Risk Assessment Methodology }
\label{sec:method}
Figure~\ref{fig:fw} denotes an overview of the proposed framework. In this section, we first discuss the structural representation of the adversarial model built upon the system models we discussed in Section \ref{sec:background} and the assumptions that affect the security of an ICS. Then, using the degradation model and its distribution (Section \ref{sec:deg}), we quantify the cyber-to-physical impact of the cyberattack. Afterward, we present the risk assessment formulation.

\subsection{Adversarial Model}
\label{sec:advmod}
To model the worst-case attacker for a defender in an ICS, we assume a rational attacker defined as follows:
\begin{definition}[Rational Attacker]
    \textit{A rational attacker intrudes into the physical system through the cyber system with limited effort and imposes physical damage to the physical system by injecting malicious data through compromised sensors and controllers using their knowledge and expertise.}
\end{definition}
This definition does not mean that we are assuming either a specific behavior of an attacker or their amount of knowledge about the system configuration. However, it helps us to study the worst-case attack scenario, which matters in designing a detection scheme.

Thinking like rational attackers, they intend to intrude on the cyber system with limited effort. To implement this intention, similar to \cite{zheng2019interdiction}, we model the cyber system as an exploitation-time attack graph (Section \ref{sec:cyber}) where each arc $(i,j)\in\A$ corresponds to a vulnerability associated with a required time $t_{ij}\in\R_{+}$ to be exploited. A rational attacker on the cyber system wants to complete an attack expeditiously, which corresponds to finding a \textit{critical path} to the target states on $\G$. Such a path is found by choosing a proper set of vulnerabilities that requires limited exploitation time and satisfies precondition dependencies. As precondition dependencies, we assume an exploitation can start if at least one of its predecessor exploitations has been completed, i.e., \textsc{OR}-type dependency. On this path, if an exploitation is successful, it can be used to reach multiple targets. Further, the attacker's capability to cause harm or disruption is constrained; thus, not all sensors and controllers of the system are under threat of simultaneous cyberattack. To capture this limitation, we consider ($i$) $T\in\Z_{+}$ as the attacker's overall available time to execute the attack actions, and ($ii$) $K\in\{1,\dots,|\N_{s}|+|\N_{c}|\}$ as the number of targets the attacker intends to compromise. We are considering a time constraint because many attackers have limited time to access different sensors and controllers, and by expanding the time constraint, we also consider a long-term attacker. The constraint on $T$ and $K$ allows us to consider different types of attackers with different levels of expertise to access the sensors and controllers.

Afterward, once a set of sensors or controllers are compromised, the attacker starts manipulating sensor measurements or control actions of the physical system. To incorporate this intention into our formulation, we utilize the state-space model (\ref{eq:nsys}) under attack that is defined as
\begin{subequations}\label{eq:advmod}
\begin{align}
    \mathbf{x}_{t+1} &= A\mathbf{x}_{t} + B[\mathbf{u}_{t}+\bb_{t}]+\mathbf{v}_{t} \\
    \mathbf{z}_{t} &=D\mathbf{x}_{t} + \ba_{t} + \mathbf{w}_{t}
\end{align}
\end{subequations}
where $\ba_{t}\in\R^{|\N_{s}|}$ and $\bb_{t}\in\R^{|\N_{c}|}$ are the attacker's manipulative actions through the compromised sensors or controllers. Since the attacker cannot perceive the exact behavior of the physical system, they alter the expected evolution of the physical system. Therefore, we ignore the Gaussian noises in (\ref{eq:advmod}). The matrices $A,B,$ and $C$ determine the attacker's knowledge of the physical system's state evolution, actuation, and measurement. Based on this knowledge, the attacker may have different attack strategies to impose various levels of physical impact.

Moreover, the state $\x_{t}$ is not directly observable to the attacker. Instead, the attacker can partially realize the sensor measurements $\z_{t}$ and determine how the physical system evolves. The attacker's restriction on realizing the sensor measurements is assumed to be
\begin{equation}\label{eq:resz}
\Vert \z_{t} \Vert_{\infty} \leq \delta.
\end{equation}
In addition, given the observed measurements $\z_{t}$, the attacker estimates the current state $\x_{t}$ and the control actions $\bu_{t}$ to decide on the attack actions. To quantify this behavior of the attacker, notice the relationship between $\z_{t}$ and $\x_{t}$ is linear in (\ref{eq:nsys}). Applying linear regression, the estimation of the states is determined as $\Hat{\x}_{t} = \left(C^{T}C\right)^{-1}C^{T}\z_{t}$. Then, the control actions $\bu_{t}$ are estimated based on $\Hat{\x}_{t}$ and $\z_{t}$. We assume the attacker estimates the control actions as a direct perception of the measurements, i.e., 
\begin{equation}\label{uez}
\bu_{t}=E \z_{t}
\end{equation}
For the attacker, $E$ should be determined such that $\bu_{t}$ preserves the system state close to some constant target value $\Tilde{\x}$. Therefore, $\bu_{t}$ is determined by minimizing $\mathbb{E}(\Hat{\x}_{t}-\Tilde{\x})$ that implies $B\bu_{t}=\Tilde{\x}-A \Hat{\x}_{t}$. Thus, assuming $B$ and $D$ have full column rank, we obtain 
\begin{equation}\label{eq:e}
    E=- B^{\dagger} A C^{\dagger}
\end{equation}
where $X^{\dagger}$ is the pseudo-inverse of $X$. In addition to the attacker's physical limitations stated above, they aim to alter the physical process so that the attack actions $\ba_{t}$ and $\bb_{t}$ bypass the existing intrusion detection mechanism. To quantify this aim, we consider a detection threshold on the attack actions $\ba_{t}$ and $\bb_{t}$, and we consider the detection algorithm as 
\begin{subequations}\label{eq:detect}
    \begin{align}
        \underline{\ba}_{i} &\leq \ba_{it} \leq \overline{\ba}_{i}, \quad i\in \N_{s} \\
        \underline{\bb}_{j} &\leq \bb_{jt} \leq \overline{\bb}_{j}, \quad j\in \N_{c}
    \end{align}
\end{subequations}

As mentioned before, gaining control over a set of sensors or controllers allows the attacker to begin manipulating the physical process. Let $\alpha_{it}\in\{0,1\}$ and $\beta_{jt}\in\{0,1\}$ represent whether the attacker has gained access to sensor $i\in \N_{s}$ and controller $j\in \N_{c}$ by time $t\in T$ or not, respectively. In other words, $\alpha_{it}$ and $\beta_{jt}$ take 1 if the target states $i\in\N_{s}$ and $j\in\N_{c}$ are compromised by time $t\in T$; otherwise, they remain 0. Once the access is gained at some time $t\in T$, the attacker begins manipulating the physical system. That means $\ba_{it}$ and $\bb_{jt}$ become nonzero for compromised sensor $i\in \N_{s}$ and controller $j\in \N_{c}$, and, otherwise, they remain zero. Accordingly, this integration is implemented by transforming the detection mechanism~(\ref{eq:detect}) into 

\begin{subequations}\label{eq:integrate}
    \begin{align}
        \underline{\ba}_{i}\alpha_{it} & \leq \ba_{it} \leq \overline{\ba}_{i}\alpha_{it}, \quad i\in \N_{s} \\
        \underline{\bb}_{j}\beta_{jt} & \leq \bb_{jt} \leq \overline{\bb}_{j}\beta_{jt}, \quad j\in \N_{c}
    \end{align}
\end{subequations}
In the proposed framework, for any $i\in\N_{s}$, $j\in\N_{c}$, and $t\in T$, $\alpha_{it}$ and $\beta_{jt}$ are considered as the decision variables that are determined by solving the framework. They represent the actions of the attacker through the cyber network; whereas $\ba_{it}$ and $\bb_{it}$ are how the attacker uses their attack to manipulate the physical process. In Section~\ref{sec:model}, we further elaborate on how the value of these variables is determined.

The most important goal of an attacker is to impose a physical impact on the physical system of an ICS by considering the cyber and physical systems and their limits. In the following sub-section, we elaborate on this objective.
\subsection{Cyber-to-Physical Impact Quantification}
\label{sec:phys}
On a physical system, the attack actions manipulate the sensor and controller data that essentially alter $\x_{t}$, accelerating the unexpected system failure. The cyber-to-physical outcome can be quantified by the unexpected failures brought on by the cyberattack. Let the random variable $\T_{f}$ denote the TTF of the system. Since the attacker perceives a realization of the system in discrete time intervals, considering the degradation model (\ref{eq:degmodel}), the TTF is defined as 
\begin{equation}
    \T_{f}=\inf\{t>0:\bs_{t}\geq\lambda\}
\end{equation}
where $\lambda$ is the known failure threshold. The attacker aims to maximize the physical impact by minimizing the MTTF, i.e., $\text{MTTF}=\E(\T_{f})$. To characterize MTTF from the attacker's viewpoint, notice that 
\begin{equation}
    \E(\T_{f}) = \sum\nolimits_{\tau=1}^{\infty} \PR(\T_{f}\geq\tau)
\end{equation}
where $\PR(\T_{f}\geq\tau)=\PR(\T_{f}\geq\tau, \bs_{\tau} \geq \lambda)+ \PR(\T_{f}\geq\tau, \bs_{\tau} < \lambda)$. From the attacker's perspective, if the system is in a failure state, attacking the system is not beneficial that implies $\{\T_{f}\geq\tau, \bs_{\tau} \geq \lambda\}=\emptyset$. Considering the attacker's restriction on the available time for executing a cyberattack on an ICS, $T$, the MTTF, for the attacker, is defined as $\E(\T_{f}) = \sum_{\tau=1}^{T} \PR(\bs_{\tau} < \lambda)$ since $\{\T_{f}\geq\tau, \bs_{\tau} < \lambda\}=\{\bs_{\tau} < \lambda\}$. 
Now, by degradation signal distribution (\ref{eq:dist}), the MTTF for the attacker can be re-written as
\begin{equation}\label{mtf}
    \E(\T_{f}) = \sum\nolimits_{\tau=1}^{T} \Phi(z_{\tau})
\end{equation}
where $\Phi(\cdot)$ is the c.d.f. of the standard normal distribution, and
\begin{equation}\label{zt}
z_{\tau} = \frac{\lambda - \sum_{t=1}^{\tau}\theta(\mathbf{x}_{t})}{\sigma_{s}\sqrt{\tau}}
\end{equation}
Since the attacker aims to minimize MTTF, by monotonicity of the c.d.f., minimizing $\E(\T_{f})$ in (\ref{mtf}) is equivalent to minimizing $\sum_{\tau=1}^{T} z_{\tau}$. We assume the degradation rate for the attacker is a direct perception of $\x_{t}$~\cite{li2021degradation}, i.e., $\theta(\x_{t})=\kappa+\gamma^{T}\x_{t}$ where $\kappa$ is constant drift, and $\gamma$ represents correlation between system states. Hence, (\ref{zt}) can be re-written as

\begin{equation}\label{zt1}
z_{\tau} = \frac{\lambda-\kappa\tau-\gamma^{T}\sum_{t=1}^{\tau}\mathbf{x}_{t}}{\sigma_{s}\sqrt{\tau}}
\end{equation}
\subsection{Worst-Case Attack Identification}
\label{sec:model}
Our proposed framework is formulated in (\ref{model}), which consists of two interconnected networks and their corresponding constraints while optimizing the attacker's objective. Table~\ref{tab:notation} lists and describes all notations used throughout the framework.

As discussed in Section \ref{sec:advmod}, the objective function (\ref{mod:obj}) incorporates the attacker's objective, which is to reach the target sensors or controllers by exploiting a limited number of vulnerabilities in the cyber system and minimize the MTTF to impose physical damage corresponding to the set of compromised sensors and controllers in the cyber system. The detailed explanation of each constraint in the framework (\ref{model}) according to their association with each cyber or physical system is as follows:

\subsubsection{Cyber System Constraints} We assume $\G$ has a single initial vulnerable state with $h_{0}=0$. This assumption is not restrictive. If multiple initial vulnerable states exist, we consider their start time zero, and we add an auxiliary node to $\G$ with $h_{aux}=0$ and add directed arcs from the auxiliary node to all initial vulnerable states with zero exploitation time.
\begin{itemize}
    \item On attack graph $\G$, since all states are of \textsc{OR}-type, in each state $i\in\N$, the attacker exploits the vulnerability that requires limited exploitation time to reach state $j$, i.e., $h_{j}=\min\nolimits_{i:(i,j)\in\A}h_{i}+t_{ij}$. This non-linear relation can be linearized by binary variables $y_{ij}$ as in (\ref{mod:cyberc1}).
    \item The attacker pays the required effort to exploit a vulnerability once but may benefit from that to reach multiple targets. This means that if the attacker does not exploit a vulnerability, no goal node is attacked through that exploitation. This limitation is modeled in (\ref{mod:cyberc3}).
    \item At each attack state $i\in\N$, the attacker hacks a set of devices on the cyber system to reach attack state $j\in\N$, such that $(i,j)\in \A$. Since the attacker has limited resources for attacking the targets, $K$, and any exploitation can be used to attack different targets, they can hack the devices based on their remaining resources at each attack state $i\in\N$. This cyber system constraint is formulated in (\ref{mod:cyberc2}). As this constraint is a network-flow constraint, we add an auxiliary node $\nu$ (super sink node) to $\G$ and add directed arcs from all target states to this node with zero exploitation time. In other words, we define the arcs set $R=\{(u,\nu):u\in\N_{s}\cup\N_{c}\}$ with $t_{ij}=0$, for $(i,j)\in R$, and $\A^{\prime}=\A\cup R$. This node can be interpreted as the attacker's imaginary target state where the physical state is under control.
\end{itemize}
\begin{subequations}\label{model}
\begin{align}
\min \ & \sum\nolimits_{(i,j)\in\A}y_{ij} \notag \\
& + \sum_{\tau=1}^{T} \frac{\lambda-\kappa\tau-\gamma^{T}\sum\limits_{t=1}^{\tau}\mathbf{x}_{t}}{\sigma_{s}\sqrt{\tau}} & \label{mod:obj}\\
\text{s.t.} \ &  h_{j} \geq h_{i}+t_{ij}-T(1-y_{ij}) &  (i,j)\in\A \label{mod:cyberc1}\\
& f_{ij}\leq K y_{ij} &  (i,j)\in\A^{\prime} \label{mod:cyberc3}\\
& 
\sum_{j:(i,j)\in\A^{\prime}}f_{ij} - \sum_{j:(j,i)\in\A^{\prime}}f_{ji} \notag \\
& \quad= 
\begin{cases}
K & i=0 \\[-3pt]
-K & i=\nu \\[-3pt]
0 & i\in\N\setminus\{0\}
\end{cases} &  i\in\N\cup\{\nu\} \label{mod:cyberc2}\\
& \alpha_{it} \leq \sum\nolimits_{j:(j,i)\in\A}f_{ji} &  i\in\N_{s},t\in T \label{mod:conc1}\\
& h_{i}-1 \leq \sum\nolimits_{t\in T}(1-\alpha_{it}) &  i\in \N_{s} \label{mod:conc2}\\
& \alpha_{i,t-1} \leq \alpha_{it} &  i\in \N_{s}, t\in T_{1} \label{mod:conc3}\\
& \beta_{it} \leq \sum\nolimits_{j:(j,i)\in\A}f_{ji} &  i\in\N_{c},t\in T \label{mod:conc4}\\
& h_{i}-1 \leq \sum\nolimits_{t\in T}(1-\beta_{it}) &  i\in \N_{c} \label{mod:conc5}\\
& \beta_{i,t-1} \leq \beta_{it} &  i\in \N_{c}, t\in T_{1} \label{mod:conc6}\\
& \x_{t+1}= A\x_{t}+B\left[\bu_{t}+\bb_{t}\right] & \qquad  t\in T \label{mod:physc1}\\
& \z_{t}= C\x_{t}+\ba_{t} & \quad  t\in T \label{mod:physc2}\\
& \bu_{t}=E\z_{t} & \quad  t\in T \label{mod:physc3}\\
& -\delta \leq \z_{it} \leq \delta &  i\in \N_{s}, t\in T \label{mod:physc4}\\
& \underline{\ba}_{i}\alpha_{it} \leq \ba_{it} \leq \overline{\ba}_{i}\alpha_{it} & \quad i\in \N_{s}, t\in T \label{mod:physc5}\\
& \underline{\bb}_{i}\beta_{it} \leq \bb_{it} \leq \overline{\bb}_{i}\beta_{it} & \quad i\in \N_{c}, t\in T \label{mod:physc6}\\
& h_{i}\in\R_{+} &  i\in \N\\
& f_{ij}\in\mathbb{Z}_{+},y_{ij}\in\{0,1\} & (i,j)\in\A^{\prime} \\
& \alpha_{it}\in\{0,1\} &  i\in \N_{s}, t\in T \\
& \beta_{it}\in\{0,1\} &  i\in \N_{c}, t\in T \\
& \x_{t}\in\mathbb{R}^{|M|},\ba_{t},\z_{t}\in\mathbb{R}^{|\N_{s}|} &  t\in T \\
& \bb_{t},\bu_{t}\in\mathbb{R}^{|\N_{c}|} & t\in T 
\end{align}
\end{subequations}

\begin{table}[!b]
\centering
\vspace{-3ex}
\caption{Notation used in the risk assessment framework.}
\label{tab:notation}
\resizebox{\linewidth}{!}{
\begin{tblr}{
  cell{1}{1} = {c=2}{},
  cell{7}{1} = {c=2}{},
  cell{21}{1} = {c=2}{},
  hline{1,32} = {-}{0.08em},
}
{\normalsize\textbf{\textsc{Sets}}}                                             &                                                                                                \\
$\N$                                             & Set of nodes (attack states) on attack graph $\G$                                              \\
$\A$                                             & {Set of directed arcs (i.e., exploitations)~on \\attack graph $\G$}                            \\
$\N_{s}\subset\N$                                & Set of nodes where sensors are compromised                                                     \\
$\N_{c}\subset\N$                                & Set of nodes where controllers are compromised                                                 \\
$M$                                              & Set of physical system states                                                                  \\
{\normalsize\textbf{\textsc{Parameters}}}                                        &                                                                                                \\
$T$                                              & Attacker's available time for attack execution                                                 \\
$K$                                              & {Number of sensors and controllers the attacker \\intends to compromise}                       \\
$t_{ij}$                                         & {Completion time of an exploit in attack state\\$i$ that its success leads to state $j$}   \\
$A\in\R^{|M|\times|M|}$                          & Attacker's knowledge of system matrix                                                          \\
$B\in\R^{|M|\times|\N_{c}|}$                     & Attacker's knowledge of input matrix                                                           \\
$C\in\R^{|\N_{s}|\times|M|}$                     & Attacker's knowledge of output matrix                                                          \\
$E\in\R^{|\N_{c}|\times|\N_{s}|}$                & Attacker's estimate of input-output matrix                                                     \\
$\gamma\in\R^{|M|}$                              & {Correlation between system states in \\degradation rate function}                             \\
$\kappa\in\R$                                    & Constant drift in degradation rate function                                                    \\
$\lambda \in\R_{+}$                              & {Attacker knowledge of failure threshold on \\degradation signal}                              \\
$\delta\in\R^{|\N_{s}|}$                         & {Attacker's limitation parameter on perceiving \\sensors data}                                 \\
$\underline{\ba},\overline{\ba}\in\R^{|\N_{s}|}$ & {Stealthiness parameters of attack actions~on \\measurements (sensors)}                        \\
$\underline{\bb},\overline{\bb}\in\R^{|\N_{c}|}$ & {Stealthiness parameters of attack actions~on \\controllers}                                   \\
{\normalsize\textbf{\textsc{Decision Variables}}}                                &                                                                                                \\
$h_{i}\in\R_{+}$                                 & {Start time of security state $i\in\N$ on\\attack graph $\G$}                                  \\
$y_{ij}\in\{0,1\}$                               & {1 if exploit $(i,j)\in\A$ is used\\to attack, o.w. 0}                                         \\
$f_{ij}\in\Z_{+}$                                & {Number of target states attacked \\through $(i,j)\in\A$}                                      \\
$\alpha_{it}\in\{0,1\}$                          & {1 if attacker obtains the full control of sensor\\$i\in \N_{s}$ by time $t\in T$, o.w. 0}     \\
$\beta_{it}\in\{0,1\}$                           & {1 if attacker obtains the full control of\\controller $i\in \N_{c}$ by time $t\in T$, o.w. 0} \\
$\x_{t}\in\R^{|M|}$                              & Physical system state at time $t\in T$                                                         \\
$\bu_{t}\in\R^{|\N_{c}|}$                        & Control actions of system at time $t\in T$                                                     \\
$\z_{t}\in\R^{|\N_{s}|}$                         & {Sensors measurements of system \\at time $t\in T$}                                            \\
$\ba_{t}\in\R^{|\N_{s}|}$                        & {Attacker's action on compromised \\sensors at time $t\in T$}                                  \\
$\bb_{t}\in\R^{|\N_{c}|}$                        & {Attacker's action on compromised \\controllers at time $t\in T$}                              
\end{tblr}
}
\end{table}

\subsubsection{Cross-layer Constraints} The constraints set (\ref{mod:conc1}-~\ref{mod:conc6}) aim at modeling the connection between the cyber and the physical systems.
\begin{itemize}
    \item The attacker can start malicious manipulations of a particular sensor or controller after a sequence of exploits in the cyber system that leads to an attack state indicating the compromise of that sensor or controller. This constraint is modeled in (\ref{mod:conc1}) and (\ref{mod:conc4}) for sensors and controllers respectively.
    \item The attacker starts the malicious actions on the physical system after the target states are activated on $\G$. This constraint is formulated in (\ref{mod:conc2}) and (\ref{mod:conc5}) for sensors and controllers, respectively.
    \item The attacker preserves the gained control on the compromised sensor or controller until the end of attack execution time $T$. This limitation is denoted in (\ref{mod:conc3}) and (\ref{mod:conc6}) for sensors and controllers, respectively. In these constraints $T_{1}=T\setminus\{0\}$.
\end{itemize}

\subsubsection{Physical System Constraints} The constraints set (\ref{mod:physc1}-~\ref{mod:physc6}) formulate the attacker's perception and restrictions on the physical system and avoiding from being detected.
\begin{itemize}
    \item The evolution of the physical system from the attacker's viewpoint (i.e., attacker's knowledge of the physical system) is captured in constraints (\ref{mod:physc1}), (\ref{mod:physc2}), and (\ref{mod:physc3}).
    \item The constraint (\ref{mod:physc4}) is the linearization of the attacker's restriction on realizing the sensor measurements.
    \item The attacker's malicious actions on the physical system need to bypass the existing cyberattack intrusion detection mechanism and remain stealthy. That is captured in (\ref{mod:physc5}) and (\ref{mod:physc6}) where $\underline{\ba},\overline{\ba}\in\R^{|\N_{s}|}$ and $\underline{\bb},\overline{\bb}\in\R^{|\N_{c}|}$ are lower and upper bounds of attacker's action on measurements and controllers, respectively. The cross-layer decision variables and constraints force the attack actions to remain zero until the corresponding sensor or controller are compromised.
\end{itemize}
Since our proposed framework is a mixed-integer linear program (MILP), powerful general solvers like Gurobi can find the optimal solution. To evaluate the effectiveness of the proposed model, the optimization model is programmed in Python using Gurobi 10.0.1 and setting its parameters, \textit{Heuristics} to 0, \textit{IntFeasTol} to $1e-09$, \textit{IntegralityFocus} to 1, and \textit{Presolve} to 0. All the experiments are conducted on a machine with CPU model 11th Gen Intel(R) Core(TM) i7-1165G7 \@ 2.80GHz and 16GB RAM, and all the experiments are solved to optimality.

\section{Computational Experiments and Results}
\label{sec:experiment}
We conduct a numerical and case study to illustrate and validate the risk assessment methodology proposed in Section \ref{sec:method}. The numerical study is intended to show the efficiency of the proposed framework. We validate our methodology on a hardware-in-the-loop (HIL) simulation test bed to demonstrate the applicability of our method in practice. For both studies, we design experiment scenarios based on the attacker's available time to execute the cyberattack $T$ and the attacker's available resources to conduct the cyberattack $K$. This way, we generate the most critical attack for attackers with different knowledge and expertise. In order to have a fair comparison of the normal system versus the system under the optimal attack, in both numerical and case studies, we simulate the physical system (\ref{eq:advmod}) for $1000$ time units given the optimal attack actions ($\ba_{t}$ and $\bb_{t}$) obtained from solving (\ref{model}) and the same Gaussian noises in the process. We provide a pseudo-code with the respective inputs in Pseudo-Code~\ref{pcode} to summarize the procedure.
\setlength{\textfloatsep}{3pt}
\begin{algorithm}[!t]
\begin{algorithmic}
\STATE{\textbf{Inputs:}}
\begin{itemize}
    \item Attacker's limitations: $T,K$ and $\delta$ in (\ref{eq:resz})
    \item Attack Graph $\G=(\N,\A)$
    \item Exploitation times $t_{ij}$, $\forall (i,j)\in\A$
    \item State-space model dynamics in (\ref{eq:nsys}): $A,B,D,\sigma_{\mathbf{v}},\sigma_{\mathbf{w}}$
    \item Degradation model parameters in (\ref{zt1}): $\kappa,\gamma,\sigma_{s}$
    \item Data control system parameters in (\ref{eq:detect}): $\underline{\ba},\overline{\ba}, \underline{\bb},\overline{\bb}$
\end{itemize}
\PROCEDURE{Attack Generation}{\textit{Inputs}}
\STATE{
\begin{enumerate}[label=\arabic*.]
    \item Compute Matrix $E$ in (\ref{eq:e})
    \item Solve model (\ref{model}) for attacker's action on both cyber and physical systems using a general purpose solver such as Gurobi
    \RETURN{Attacker's optimal actions: 
    \begin{itemize}
        \item [ ] Cyber system: $y^{*}_{ij},f^{*}_{ij}$
        \item [ ] Physical system: $\ba^{*}_{t},\bb^{*}_{t}$
    \end{itemize}}
\end{enumerate}
}
\ENDPROCEDURE
\PROCEDURE{Risk Assessment}{$y^{*}_{ij},f^{*}_{ij},\ba^{*}_{t},\bb^{*}_{t}$}
\STATE{
\begin{enumerate}[label=\arabic*.]
    \item Deploy the attacker's optimal actions $y^{*}_{ij},f^{*}_{ij},\ba^{*}_{t},\bb^{*}_{t}$ into the ICS to identify the critical vulnerabilities in the cyber system and evaluate MFT of the physical system using (\ref{mtf}) and (\ref{zt1})
\end{enumerate}
}
\ENDPROCEDURE
\end{algorithmic}
\caption{Risk Assessment Procedure}\label{pcode}
\end{algorithm}

\begin{figure*}[!t]
    \centering
    \includegraphics[width=.8\textwidth]{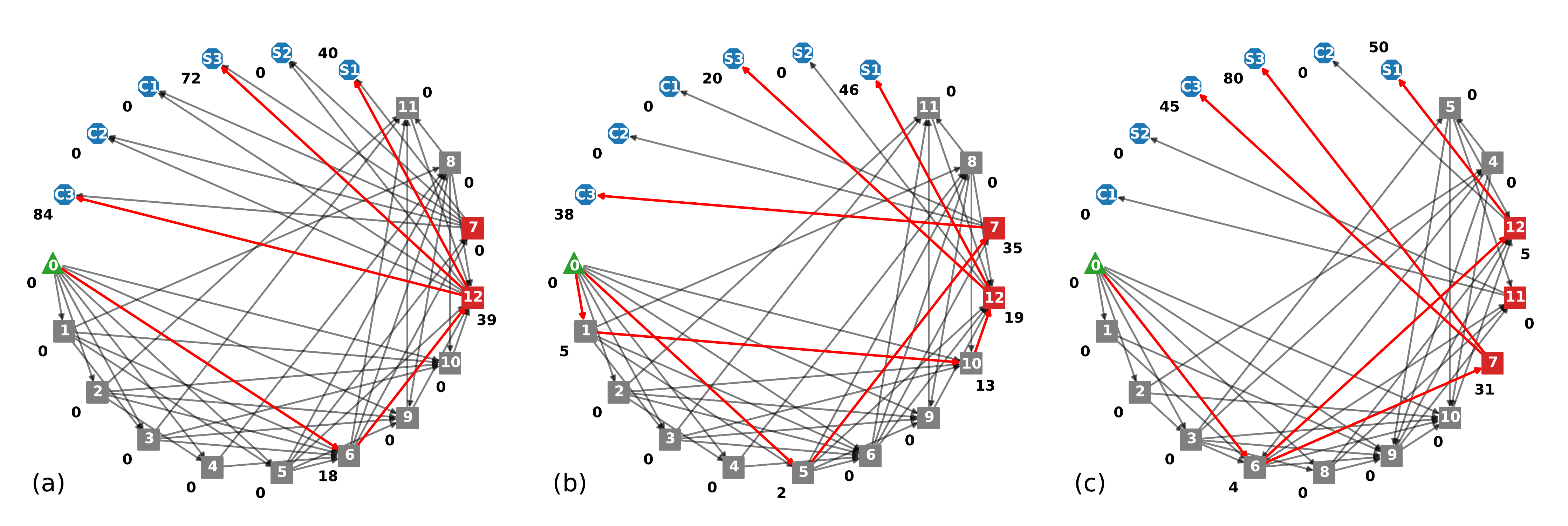}
    \caption{The generated random attack graph topologies are based on the degree of vulnerabilities. For $K=3$, the red edges represent worst-case attack actions on the cyber system, and the numbers next to each node are the time the security state is available to the attacker.}
    \label{fig:topo}
\end{figure*}

\subsection{Numerical Study}
\subsubsection{Experiment Setup}
\label{sec:setuprand}
We consider an ICS where a physical system with $|M|=3$, $|\N_{s}|=3$, and $|\N_{c}|=3$ (2nd scenario in \cite{li2021degradation}). As the cyber network, a communication network is connected to these sensors and controllers. Random cyber networks are generated with a single initial vulnerability point, 12 security states, three states for sensors and controllers representing each of them are being compromised, and $t_{ij}\sim U(1,50)$. Three topologies of random networks based on degree of vulnerabilities (the attack graph connectivity) are generated, ($i$) high degree of vulnerabilities: all sensors' and controllers' states are connected to the same precondition states (Figure~\ref{fig:topo}a), ($ii$) medium degree of vulnerabilities: sensors' and controllers' states are connected to two separate precondition states (Figure~\ref{fig:topo}b), and, ($iii$) low degree of vulnerabilities: three random combinations of sensor-controller states are connected to three separate precondition states (Figure~\ref{fig:topo}c). The physical system dynamics is given by, $B=C=\mathbf{I}_{3}$, $\sigma_{\mathbf{v}}=0.1$, $\sigma_{\mathbf{w}}=\sqrt{0.001}$, $\mathbf{v}_{t}\sim\mathcal{N}(0,\sigma_{\mathbf{v}}^{2}\mathbf{I}_{3})$, $\mathbf{w}_{t}\sim\mathcal{N}(0,\sigma_{\mathbf{w}}^{2}\mathbf{I}_{3})$, $\kappa=0.4713\sqrt{2}$, $\gamma=\begin{bmatrix} 0.058 & 0.058 & 0.996 \end{bmatrix}^{T}$, $\sigma_{s}=0.1$, $E$ is computed by the relation (\ref{eq:e}), and $$
A=
\begin{bmatrix} 0.12 & 0.30 & 0.06 \\ 0.06 & 0.48 & 0.06 \\ 0.30 & 0.12 & 0.54 \end{bmatrix}
$$
Given these data, we simulate the normal system for 50 replications with $1000$ time units and compute $\bs_{t}$. We set $T_{\lambda}=600$ to be the time the system fails under normal operation and obtain $\lambda$ by taking the average of the degradation signals at this time. Also, we compute the average standard deviation of the sensor reading over the replications, $\sigma_{\z}$ and set $\delta=3\sigma_{\z}$, $\underline{\ba}_{i}=-3\sigma_{\mathbf{w}}$, $\overline{\ba}_{i}=3\sigma_{\mathbf{w}}$, $\underline{\bb}_{i}=-3\sigma_{\mathbf{v}}$, and $\overline{\bb}_{i}=3\sigma_{\mathbf{v}}$. Given these settings, we design scenarios based on the attacker's available ($i$) resources for exploiting vulnerability $(K)$, and ($ii$) time on executing attacks $(T)$.

\begin{figure*}[!t]
    \centering
    \begin{subfigure}[!t]{0.34\linewidth}
    \includegraphics[width=\linewidth]{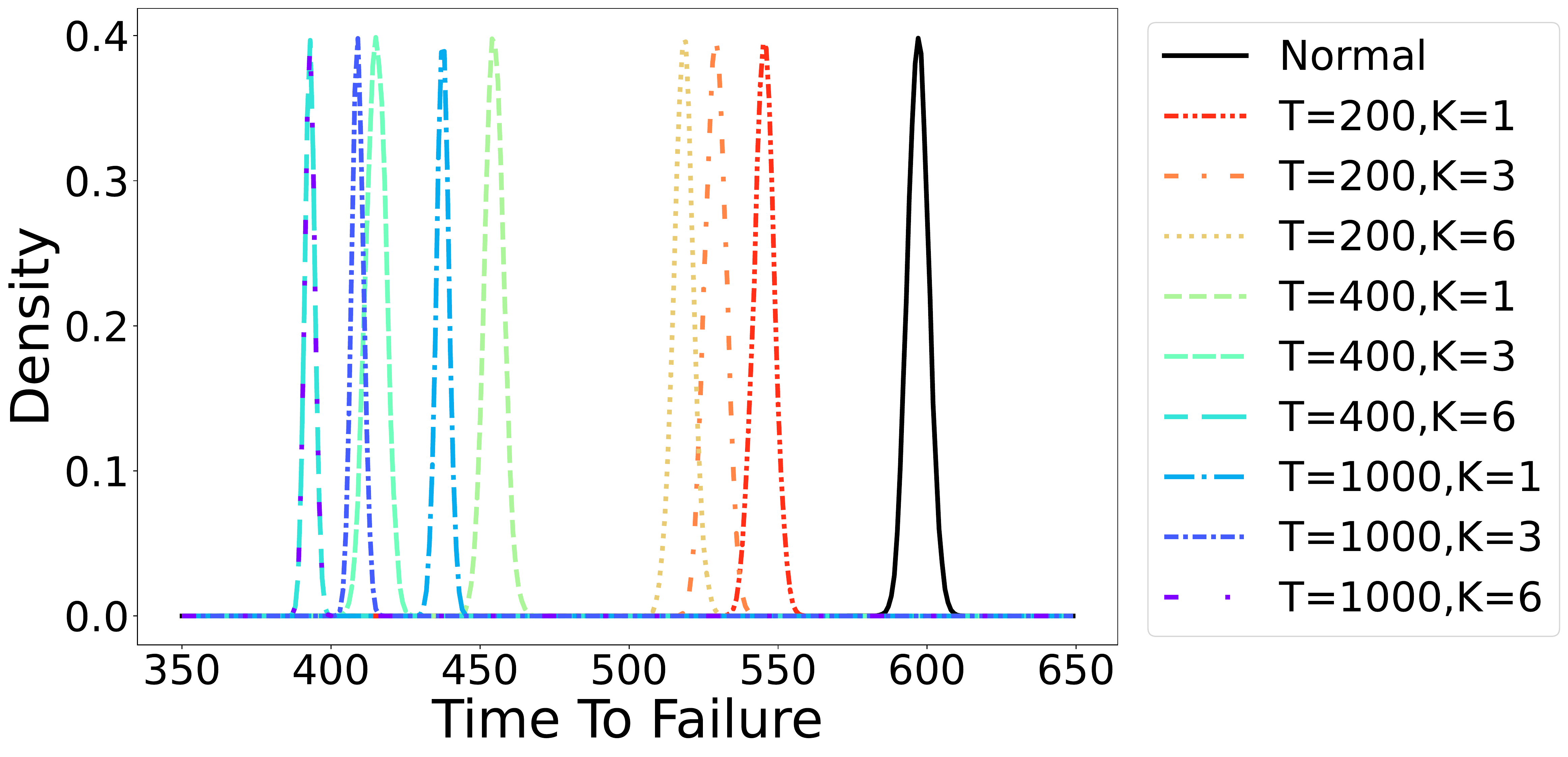}%
    \caption{Numerical study}
    \label{fig:ft_ns}
    \end{subfigure}%
    \begin{subfigure}[!t]{.34\linewidth}
    \includegraphics[width=\linewidth]{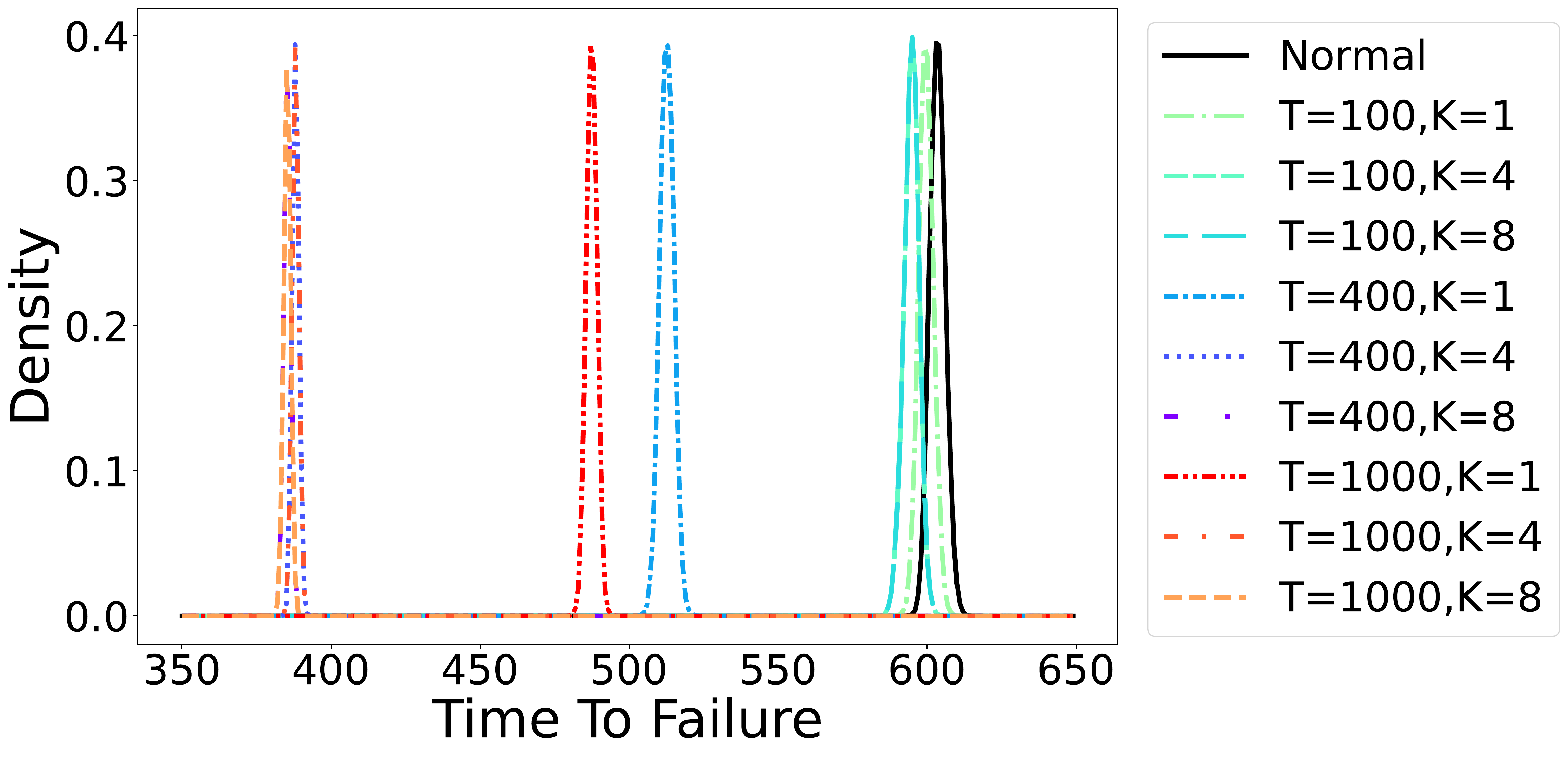}
    \caption{Case study}
    \label{fig:ft_case}
    \end{subfigure}%
    \begin{subfigure}[!t]{.32\linewidth}
        \vspace{2ex}
        \includegraphics[width=\linewidth]{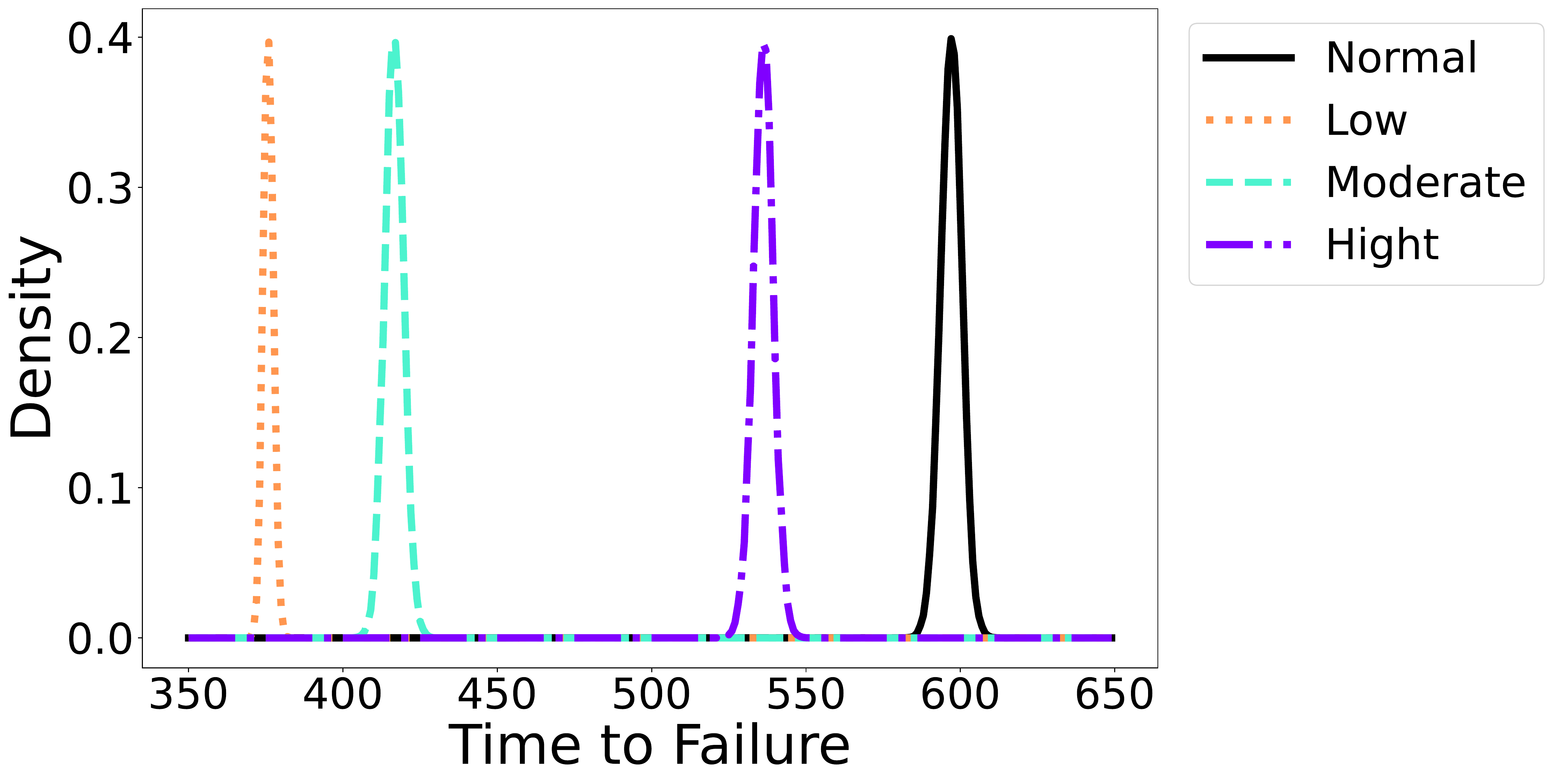}
        \caption{Sensitivity analysis on attack stealthiness when $T=400$ and $K=3$.}
        \label{fig:sa_pdf}
    \end{subfigure}
    \caption{Time to failure $\T_{f}$ distribution.}
    \vspace{-3ex}
\end{figure*}

\begin{figure}[!b]
    \centering
    \begin{subfigure}[!t]{.5\linewidth}
    \includegraphics[width=\linewidth]{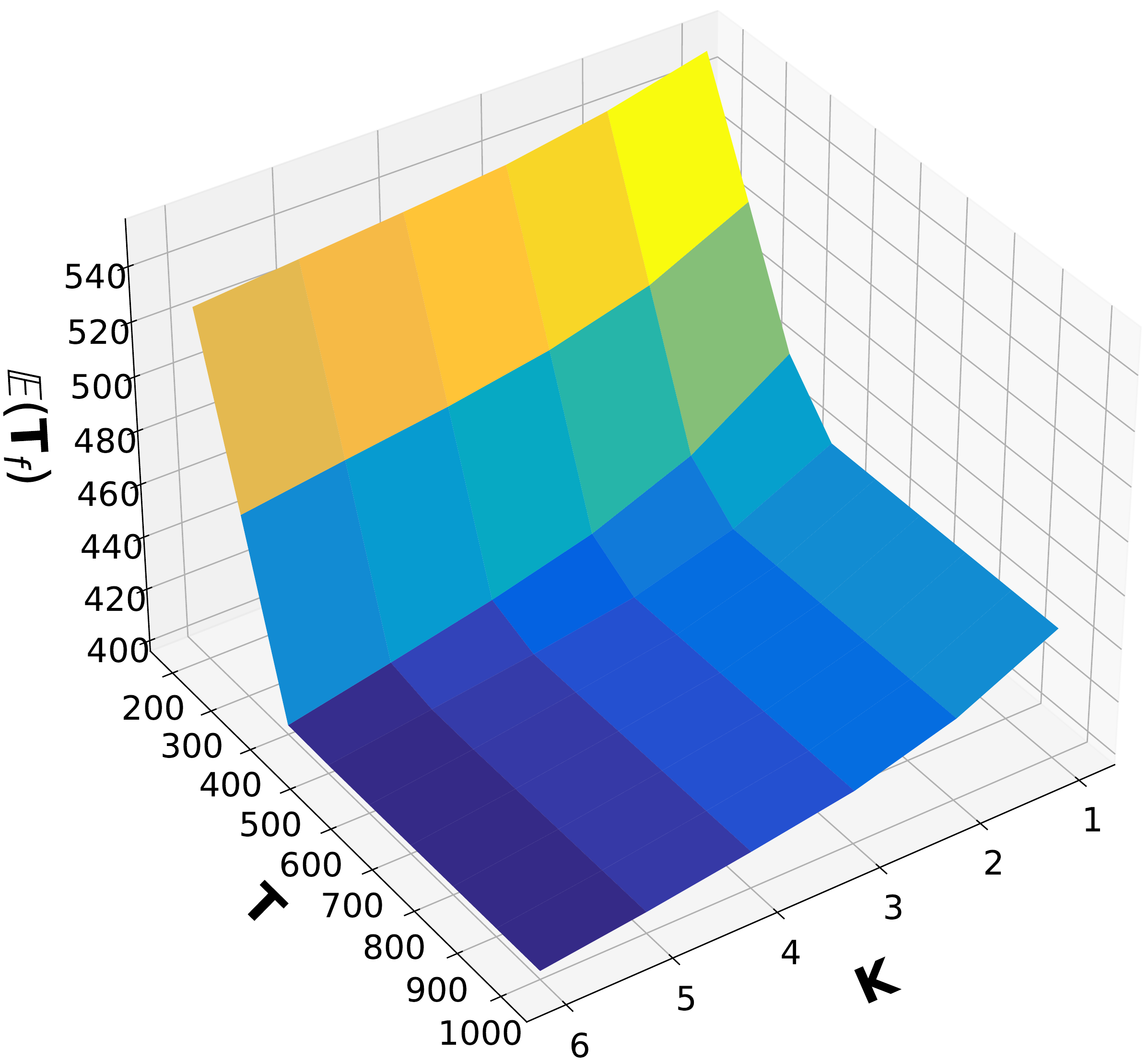}%
    \caption{Numerical study}
    \end{subfigure}%
    \begin{subfigure}[!t]{.5\linewidth}
    \includegraphics[width=\linewidth]{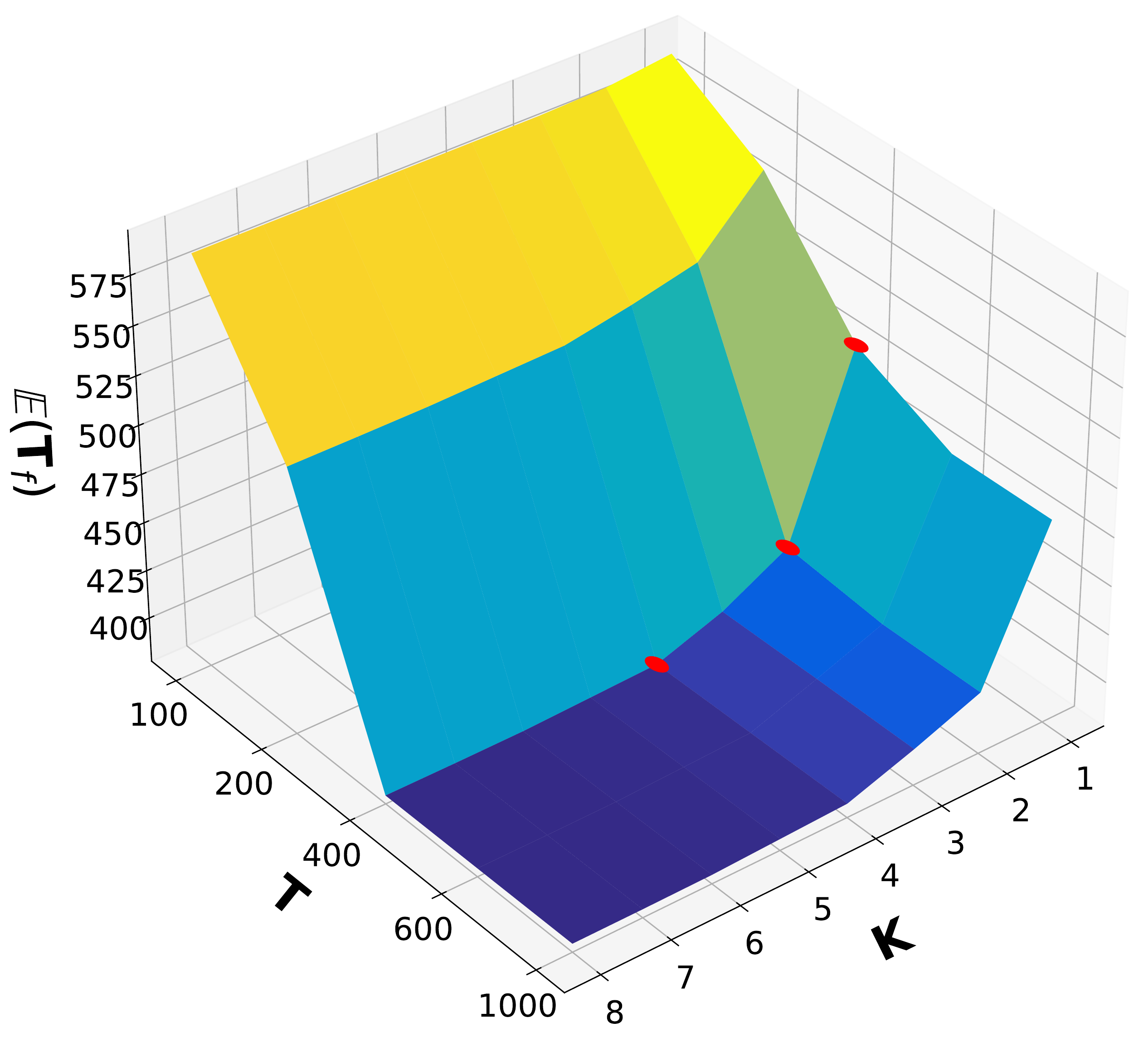}
    \caption{Case study}
    \end{subfigure}
    \caption{Mean time to failure $\E(\T_{f})$.}
    \label{fig:MTTF}
\end{figure}
\begin{figure}
    \centering
    \includegraphics[width=.6\linewidth]{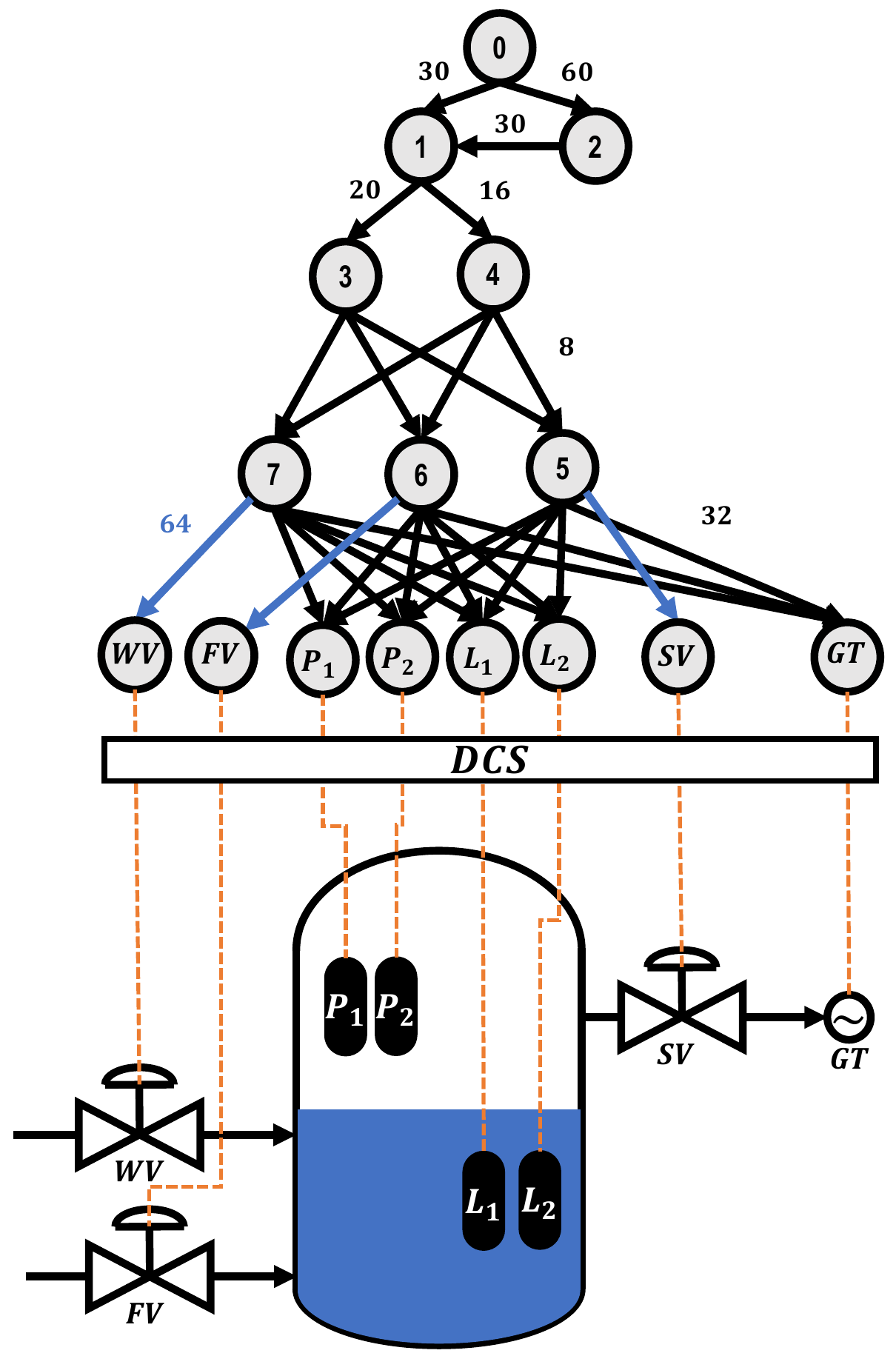}
    \caption{Cyber-physical system for the case study. Blue edges on the attack graph have the same exploitation time.}
    \label{fig:casestruc}
\end{figure}

\begin{figure*}[!t]
    \centering
    \includegraphics[width=.75\linewidth]{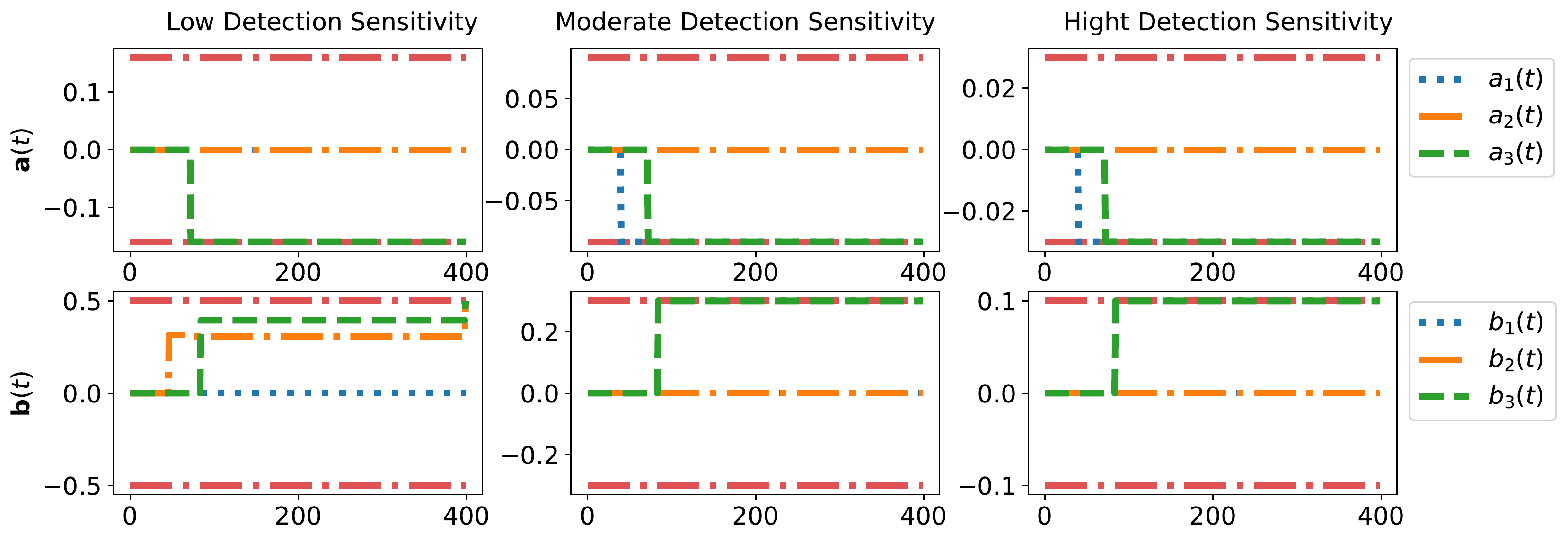}
    \caption{Worst-case attack actions on numerical study for sensitivity analysis on attack stealthiness when $T=400$ and $K=3$.}
    \label{fig:sa_act}
\end{figure*}

\begin{figure*}[!t]
    \centering
    \includegraphics[width=.75\linewidth]{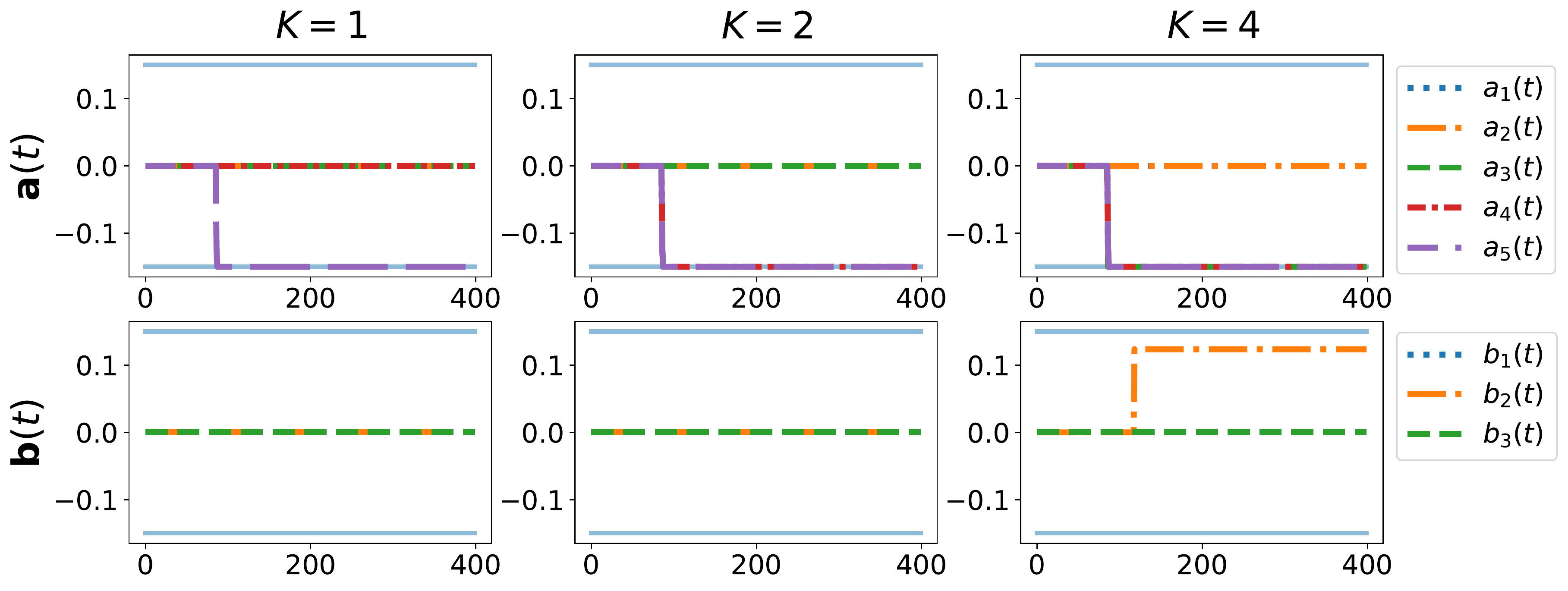}
    \caption{Worst-case attack actions on sensors and controllers of the BWPP for $T=400$.}
    \label{fig:caseact}
\end{figure*}

\subsubsection{Results}
\label{sec:rndresult}
For each random attack graph topology, the worst-case attack actions on the cyber system and the attack state start times, $h_{i}$, are denoted in Figure~\ref{fig:topo} when $K=3$. In this figure, the red edges denote the critical vulnerabilities in the cyber system whose exploitation leads to the highest decrease of MTTF when the attacker has limited resources, i.e., $K=3$. Malicious data injection starts when an attacker compromises sensors or controllers. Considering attack graph in Figure~\ref{fig:topo}, the attacker compromises the controller $C_{3}$ at $t=84$; Therefore, $\bb_{C_{3},t}=0$, for $t<84$, and, after that, $\bb_{C_{3},t}$ is determined by the attacker. Despite existing sensors and controllers that were easier to reach regarding exploitation time, the attacker targeted those that led to a noticeable physical impact. This observation denotes that, from the attacker's perspective, an ICS's cyber and physical systems are not independent, and the attacker's decision to exploit the vulnerabilities in the cyber system is correlated with the physical impact that will be imposed on the physical system. Furthermore, from Figure \ref{fig:topo}, we observe that the different topology of the attack graph leads to different times that the attacker can reach the physical system. Hence, implementing the mitigations in the cyber system should not only focus on increasing the difficulty of reaching the most accessible targets, but it also is essential to focus on the impacts that a cyberattack can impose on the physical system.

Figure~\ref{fig:ft_ns} denotes the distribution of the TTF for the physical system under normal operation and the system under worst-case attack when $K=1,3,6$ and $T=200,400,1000$. We observe that the failure rate increases as $K$ and $T$ increase, which shows how the worst-case malicious actions accelerate the TTF. However, the increase in $K$ has more effect on increasing the failure rate. Figure~\ref{fig:MTTF}(a) elaborates further on this observation. Regardless of the attacker's available time $T$, as $K$ increases, the failure rate increases and, consequently, the attacker's MTTF decreases. However, when $T\geq500$, the TTF acceleration remains almost the same as $K$ increases. This observation denotes that if the attacker has limitations on compromising the sensors and controllers but can execute a cyberattack for a longer time, the physical consequences will not be significant, and, accordingly, the chance of being detected by the intrusion detection system will increase.

\subsubsection{Sensitivity Analysis on Attack Stealthiness}
\begin{table}[!thb]
\centering
\caption{Attack stealthiness scenarios.}
\label{tab:steal}
\begin{tblr}{
  cells = {c},
  cell{1}{1} = {r=2}{},
  cell{1}{2} = {c=4}{},
  vline{2} = {1-5}{},
  hline{1,6} = {-}{0.08em},
  hline{2} = {2-5}{},
  hline{3} = {-}{},
}
{  Detection Sensitivity\\Level} & Bounds                  &                        &                         &                        \\
                                 & $\underline{\ba}$       & $\overline{\ba}$       & $\underline{\bb}$       & $\overline{\bb}$       \\
Low                              & $-5\sigma_{\mathbf{w}}$ & $5\sigma_{\mathbf{w}}$ & $-5\sigma_{\mathbf{v}}$ & $5\sigma_{\mathbf{v}}$ \\
Moderate                         & $-3\sigma_{\mathbf{w}}$ & $3\sigma_{\mathbf{w}}$ & $-3\sigma_{\mathbf{v}}$ & $3\sigma_{\mathbf{v}}$ \\
High                             & $-\sigma_{\mathbf{w}}$  & $\sigma_{\mathbf{w}}$  & $-\sigma_{\mathbf{v}}$  & $\sigma_{\mathbf{v}}$  
\end{tblr}
\end{table}

Different ICSs have various detection strategies with diverse accuracy levels. This subsection aims to analyze the sensitivity of the generated worst-case attack actions on the detection power of the ICS. This analysis provides insight into how a low detection sensitivity level provides degrees of freedom and a high detection sensitivity level restricts the attacker in choosing the attack actions. The accuracy and power of the detection mechanism are captured by the stealthiness constraints~(\ref{mod:physc5}) and~(\ref{mod:physc6}) in our framework. Considering the attack graph topology (a) in Figure~\ref{fig:topo} and the same dynamics of the physical system (see Section~\ref{sec:setuprand}), we consider three detection mechanisms, as shown in Table~\ref{tab:steal}, and solve~(\ref{model}) for worst-case attacks under each of these scenarios that determines the lower and upper bounds of the stealthiness constraints~(\ref{mod:physc5}) and~(\ref{mod:physc6}). Figure~\ref{fig:sa_pdf} depicts the TTF distribution for the physical system under normal operation and the system under attack, considering each detection sensitivity level. We observe that as the sensitivity of the detection mechanism increases, the worst-case attack actions cannot increase the failure rate and accelerate the TTF significantly. Further, we observe a noticeable increase in the failure rate of the physical system as detection sensitivity decreases. Figure~\ref{fig:sa_act} elaborates further on this observation by magnifying the worst-case attack actions. We observe that a low detection sensitivity level gives more degrees of freedom to the attacker to choose different combinations of sensors and controllers to compromise even though they have the same available knowledge and resources. In addition, with a moderately or highly sensitive detection mechanism established, the same sensors and controllers are compromised. However, as the attacker is less restricted in choosing the malicious data with a moderate detection sensitivity, the failure rate is increased noticeably. This result does not imply that a higher detection sensitivity necessarily protects the system. This is because as the detection sensitivity level increases, the false alarm rate increases when there is no attack. Therefore, a powerful and accurate detection mechanism in the physical system is required.

\subsection{Case Study}

\subsubsection{Experiment Setup}
The HIL case study described in \cite{huang2018assessing} is investigated here as an ICS comprising both cyber and physical systems. As the physical system, a boiling water power plant (BWPP) contains three states to be controlled (drum pressure, electric output, and fluid density), three controllers, a water valve (WV), a fuel valve (FV), a steam valve  (SV), and five sensors, two pressure sensors (P1 and P2), two water level sensors (L1 and L2), and a sensor (GT) for reading the generated electricity in the field area. The cyber system consists of an administrator host (AH) and two PCs in a corporate network, a human-machine interface (HM) and a data server (DS) in a supervisory network, and three embedded controllers in a control network. Given the vulnerability list and their exploitation times for this case, we construct the corresponding attack graph of the cyber system.

\begin{figure*}[!t]
    \centering
    \includegraphics[width=\linewidth]{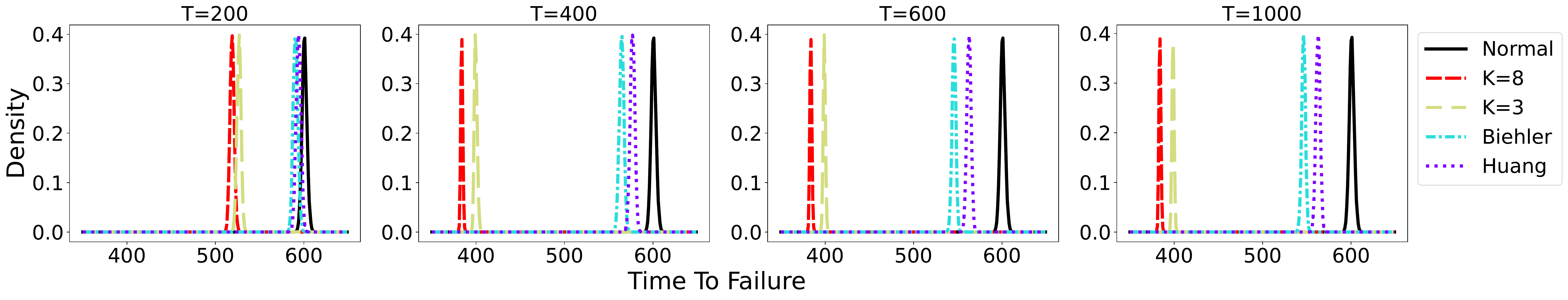}
    \caption{TTF distribution of attack actions: Huang et al.~\cite{huang2018assessing} attack scenario versus worst-case attacks obtained by our framework ($K=3$); Biehler et al.~\cite{biehler2023sage} attack generation model versus worst-case attacks obtained by our framework ($K=8$).}
    \label{fig:comppdf}
\end{figure*}

\begin{figure*}[!t]
    \centering
    \includegraphics[width=\linewidth]{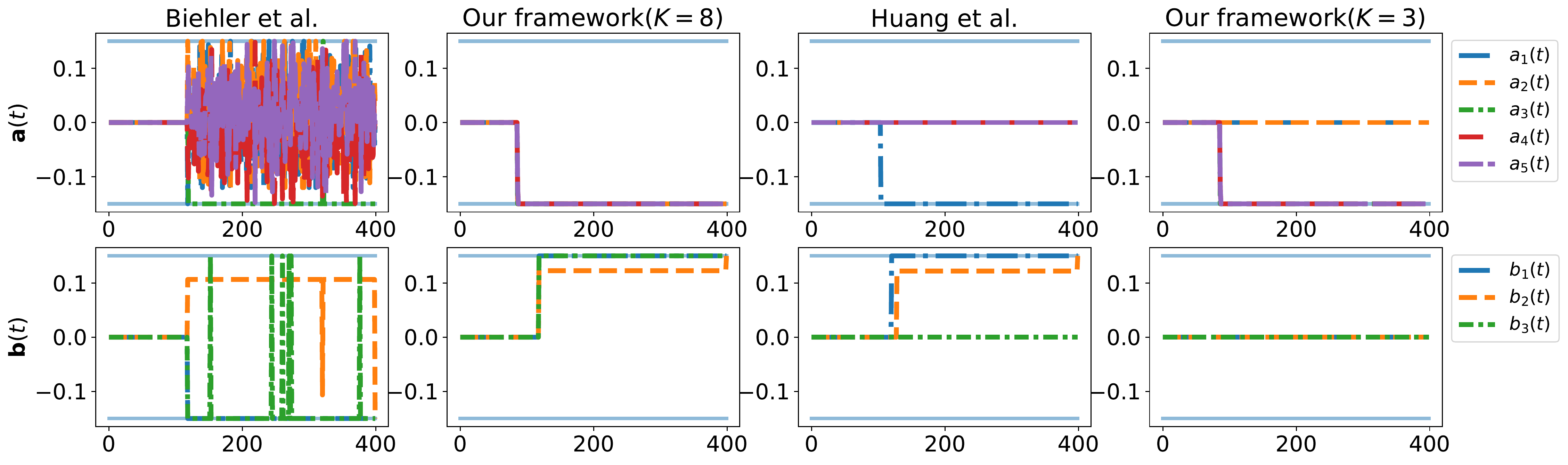}
    \caption{Attack actions on sensors and controllers: Huang et al.~\cite{huang2018assessing} attack scenario versus worst-case attacks obtained by our framework ($K=3$); Biehler et al.~\cite{biehler2023sage} attack generation model versus worst-case attacks obtained by our framework ($K=8$).}
    \label{fig:compact}
\end{figure*}

Figure~\ref{fig:casestruc} denotes the cyber-physical structure of this case. The graph represents different states, starting from initial access (`0') to acquiring user privilege on AH and PCs (`1' and `2'), root privilege on HM and DS (`3' and `4'), and successfully acting as HM and sending forged Modbus commands (`5', `6', and `7'). The remaining states correspond to the compromised sensors and controllers. Huang et al. \cite{huang2018assessing} reported the same vulnerability and exploitation time for sensors and controllers. However, compromising controllers requires more effort, so we increase their exploitation time by a factor of two. The shortest path length to compromised sensors is 86 time units, while for controllers, it is 118.

Given the physical system dynamics, we define $\Tilde{\x}_{t}=\x_{t}-\overline{\x}$ where $\overline{\x}$ is the reported steady state of the system. Since $\Tilde{\x}_{t}$ exists, and it must be that $\E(\Tilde{\x}_{t})=0$, any linear function $f(\Tilde{\x}_{t})=\beta_{0}+\beta^{T}\Tilde{\x}_{t}$ exists such that $\E(f(\Tilde{\x}_{t}))=\beta_{0}$. As discussed in Section ~\ref{sec:deg}, the degradation rate is a function of $\x_{t}$. Therefore, we assume that the exact values of $\kappa$ and $\gamma$ of this case depend on the plant type, and further investigations are out of the scope of this study. Hence, we generate random values between 0 and 1 as degradation parameters for simplicity and without losing generality. Also, we consider $\sigma_{s}^{2}=0.01$. In the same way as the numerical study, we obtain $\lambda$, $\delta$, $\underline{\ba}_{i}$, $\overline{\ba}_{i}$, $\underline{\bb}_{i}$, and $\overline{\bb}_{i}$. It is worth mentioning that we do not conduct physical experiments, but we use the dynamics provided in~\cite{huang2018assessing} in our risk assessment model in (\ref{model}).

\subsubsection{Results}
\label{sec:caseresult}
Figure~\ref{fig:ft_case} denotes the TTF distribution for the normal system operation and the system under worst-case attack when $K=1,4,8$ and $T=100,400,1000$. We observe the same results as the numerical study, and it validates our observations on a real system. However, when $T=100$, the failure rate distribution is close to the normal system regardless of the increase in $K$. Since the shortest path length to the controllers is higher than the attacker's available time, the attacker could only compromise the sensors (depending on $K$). This observation denotes the importance of protecting controllers from cyberattacks. Even though when the attacker has high resources ($K\geq6$), the controllers remain secure, and consequently, the imposed physical impact is insignificant. 
        
Figure~\ref{fig:MTTF}(b) elaborates more on this observation. We observe a noticeable jump in MTTF when $T$ increases from 200 to 400 as $K$ increases. Figure~\ref{fig:caseact} magnifies the worst-case attack actions that led to this observation. We observe that as $K$ increases, the attacker could compromise the controller $SV$, leading to aggressive actions on the sensors $GT,L_{1},$ and $L_{2}$. Furthermore, from $K=2$ to $K=4$, the attacker causes significant reduction in MTTF by compromising $SV$ and $GT$ and manipulating $L_{1}$ and $L_{2}$ in the same way as $K=2$.

\begin{figure}[!t]
    \centering
    \includegraphics[width=.5\linewidth]{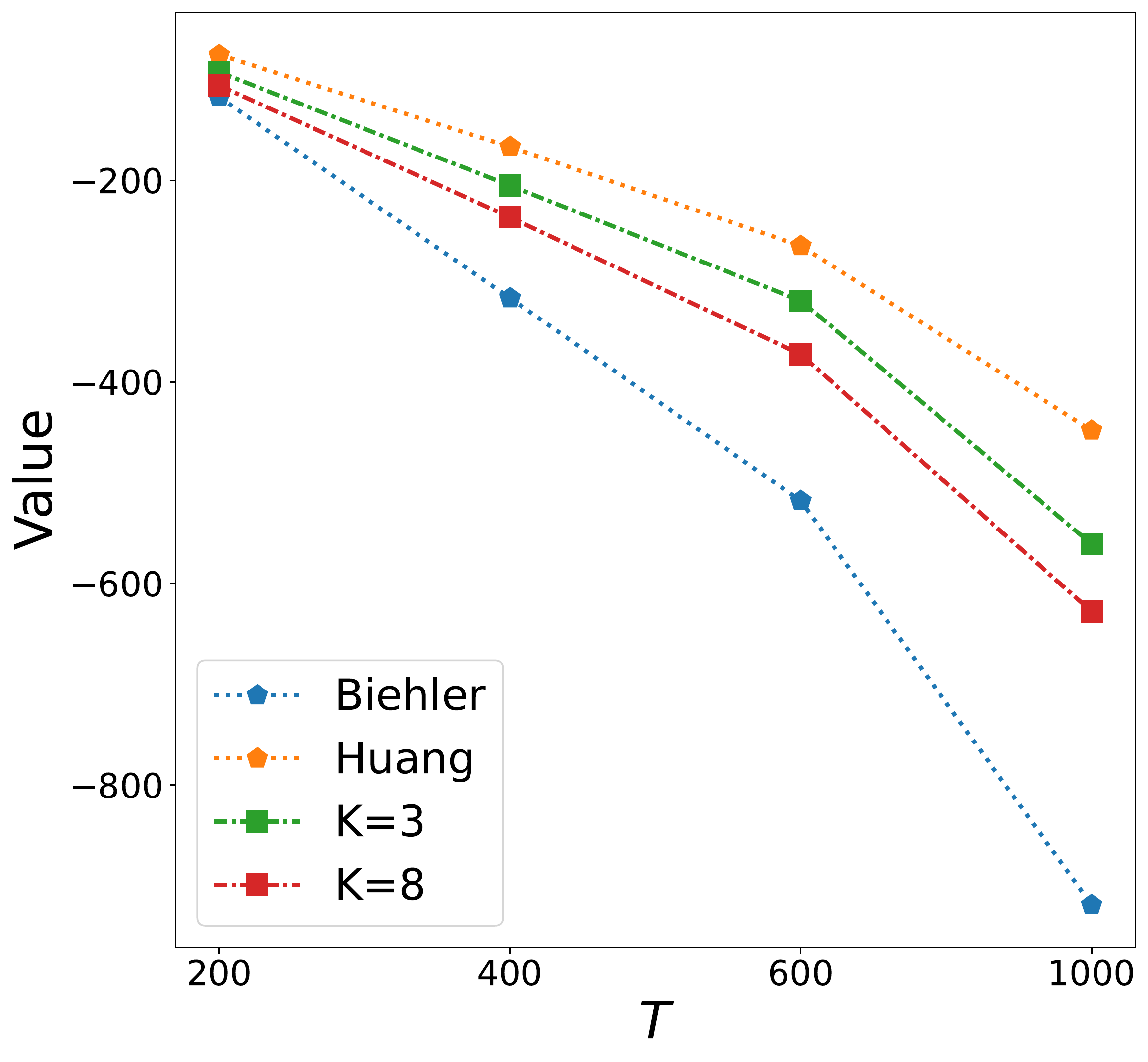}%
    \includegraphics[width=.49\linewidth]{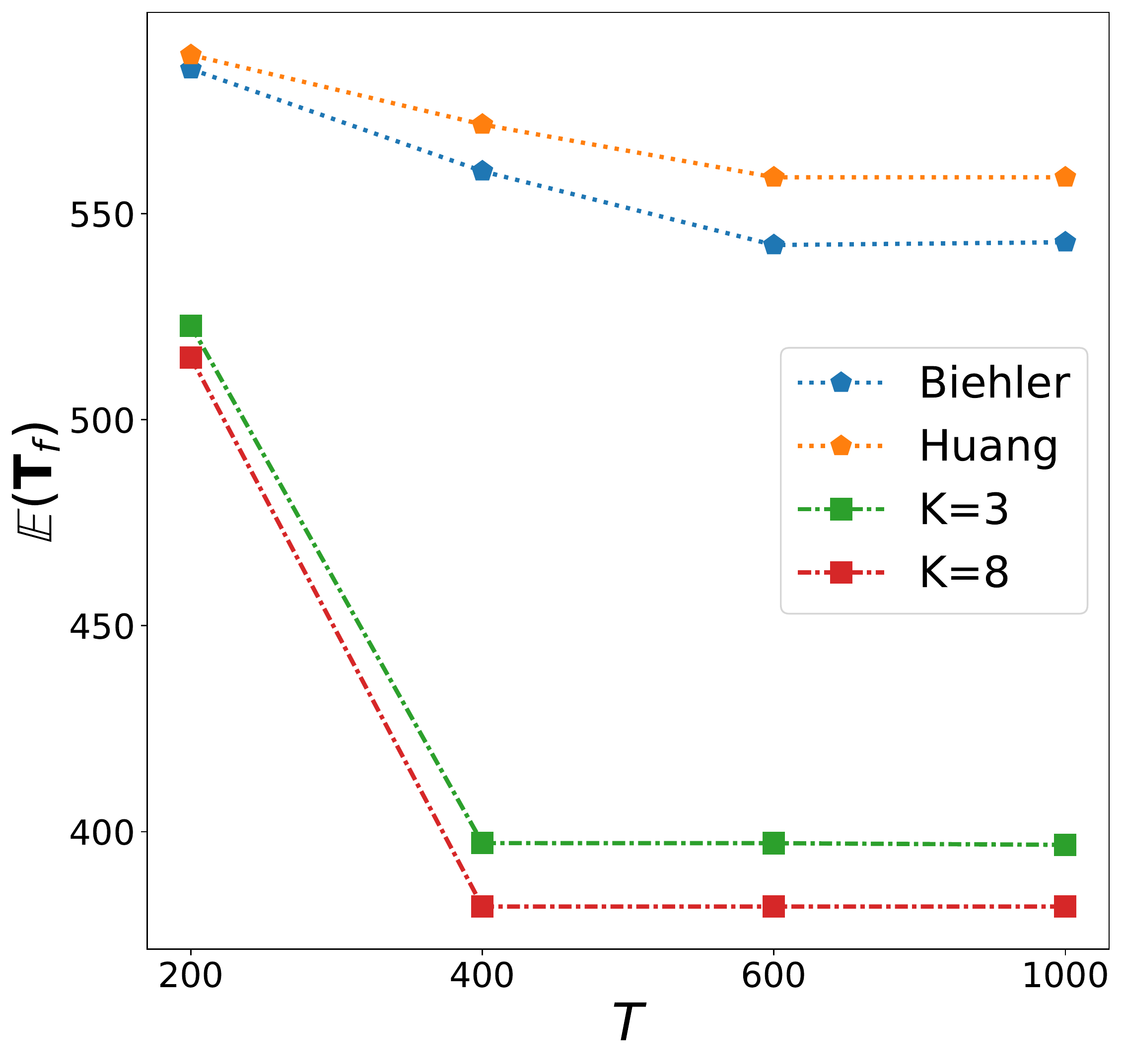}
    \caption{Damage function~(\ref{objb}) value (left) and $E(\T_{f})$ (right) of attack actions: Huang et al.~\cite{huang2018assessing} attack scenario versus worst-case attack obtained by our framework ($K=3$); Biehler et al.~\cite{biehler2023sage} attack generation model versus worst-case attacks obtained by our framework ($K=8$).}
    \label{fig:compmtf}
\end{figure}

\subsubsection{Comparative Analysis}\label{sec:comp}
Due to the different nature of the cyber-risk assessment methods in the literature, it is hard to find an ICS containing the required cyber and physical system information to establish a fair comparison. Therefore, our intention in this section is not to show that our framework outperforms the other methods. Instead, we point out the significant characteristics of an ICS that are incorporated in our formulation.

Huang et al.~\cite{huang2018assessing} considered an attack scenario to validate their method. In this scenario, the attacker follows the path $0\to1\to3\to(5,6)\to(P_{1},SV,FV)$ on the attack graph (Figure~\ref{fig:casestruc}) and compromises $P_{1},FV,$ and $SV$ at time 104, 120, 128, respectively. Then, the attacker starts attack actions on the BWPP through $P_{1},FV,$ and $SV$. This scenario is equivalent to a case where $K=3$ for some attack execution time $T$. For $T\in\{200,400,600,100\}$, we implement this attack scenario on our framework in~(\ref{model}). Also, we solve (\ref{model}) for the worst-case attack when $K=3$. The worst-case attack found by our framework follows the path $0\to1\to4\to5\to(L_{1},L_{2},GT)$ on the attack graph (Figure~\ref{fig:casestruc}) and compromises $L_{1},L_{2},$ and $GT$ at time 86. Figures~\ref{fig:comppdf} and~\ref{fig:compmtf} denote the TTF distribution and the MTTF of BWPP under Huang's attack actions versus the worst-case attack obtained from our framework when $K=3$. We observe that the MTTF of the BWPP is further decreased through the attack actions obtained from our framework. Figure~\ref{fig:compact} magnifies this observation by looking at the attack actions. We observe that the worst-case attack with $K=3$ can significantly reduce MTTF by only compromising the sensors $L_{1},L_{2},GT$ that require less effort than controllers. These results verify the inability of the probabilistic risk assessment methods to systematically identify the worst-case attack.

Bielher et al.~\cite{biehler2023sage} proposed a generic attack generation framework on the physical system without considering the cyber system and assuming the attacker has full access to the attack channels, which is unrealistic. Their approach is equivalent to the case where $K=8$ and some value of $T$ in our formulation. Our physical constraints (\ref{mod:physc1})-(\ref{mod:physc6}) correspond to the constraints in their model with special case functions (i.e., linear functions) and no cost on physical attack actions. To elaborate more on this, we formulate the BWPP as (\ref{eq:nsys}) that can be re-written as $\z_{t}=CB\bu_{t}+CA\x_{t}+\mathbf{v}^{\prime}_{t}$ where $\z_{t}$ is considered as the output of this process. By (\mbox{\ref{eq:advmod}}), the BWPP process output under attack is $\z_{t}^{A}=\z_{t}+ CB\bb_{t} + \ba_{t} +\mathbf{v}^{\prime}_{t}$. Therefore, the equivalent objective function (damage function) in their model is
\begin{equation}\label{objb}
    \min_{\ba_{t},\bb_{t}} \ -\sum_{t=1}^{T}\left\Vert \z_{t}^{ref} - CB\bb_{t} - \ba_{t} \right\Vert_{1}
\end{equation}
where $\z_{t}^{ref}$ is the sensor measurements under normal operations. In their model, the stealthiness constraint is equivalent to (\ref{eq:detect}), and the physical limits constraint corresponds to (\ref{eq:resz}) and (\ref{uez}) in our framework.

Noting the equivalence between our physical system constraints and the attack generation model in \cite{biehler2023sage}, we solve our framework (\ref{model}) for $K=8$ and $T\in\{200,400,600,1000\}$ by replacing (\ref{objb}) with (\ref{mod:obj}) to generate Biehler's attack actions. Figures~\ref{fig:comppdf} and~\ref{fig:compmtf} denote the TTF distribution and the MTTF of BWPP under Biehler's attack actions versus the worst-case attack obtained from our framework when $K=8$. We observe that our worst-case attack actions further reduce the MTTF of the BWPP. The worst-case attacker with limited knowledge and resources (i.e., $K=3$) reduce the MTTF of significantly more than the attack actions generated by Biehler's model. Figure~\ref{fig:compact} elaborates more on this observation. The attack actions identified by our framework noticeably differ in pattern from those generated by Biehler's approach. This difference arises from using varied objectives to assess the physical consequences of the attack actions. Figure~\ref{fig:compmtf} (left) denotes the damage function value (\ref{objb}) in Biehlers's formulation for the worst-case attack actions, for $K=3$ and $K=8$, attack actions obtained from Biehlers's model, and Huang's attack scenario. The objective function (\ref{objb}) essentially identifies attack actions that result in outputs with the highest $l_{1}-$norm distance from the reference sensors measurement. We observe that Biehler's attack actions generate attacks that alter the measurements more than the other approaches, even though these actions do not significantly affect MTTF.

The results of this comparative analysis highlight our framework's unique advantages and practical limitations. The primary advantages include: ($i$) Our framework diverges from the probabilistic risk assessment methods, such as~\cite{huang2018assessing}, which only quantifies the physical impact of a given attack, by modeling worst-case attacks under the attacker's rationale. This approach not only identifies the highest-impact attacks among all possible attacks but also provides a holistic view of the system's weaknesses. This comprehensive understanding of risks, under limitations on the attacker's knowledge and resources, equips operators with a deeper insight into the resilience of their ICS. ($ii$) By jointly optimizing attack actions across both cyber and physical layers, our method advances beyond existing models by identifying the critical components of the ICS. This holistic characterization of attack actions enables us to surpass methods that only consider the physical layer, as in~\cite{biehler2023sage}, by eliminating the assumptions on the accessibility of all attack channels. Such detailed identification is pivotal for operators to strategize the design of detection and mitigation schemes. Nonetheless, there are certain limitations to our approach. Our framework does not incorporate uncertainties in the success rate of cyber exploitation and the stochastic nature of physical processes. This aspect is crucial for a more realistic risk assessment and is a focus for our future research. Another limitation is that the effectiveness of our model is contingent on having detailed knowledge of cyber vulnerabilities (to build the attack graph) and physical process parameters (to implement the constraints and model the degradation). This requirement can pose a challenge in cases where such information is limited.

\section{Conclusion}
\label{sec:conclusion}
In this article, we introduce an integrated framework for identifying the worst-case cyberattacks on ICSs and assess their physical consequences, combining cyber and physical system. We model the attacker's strategy as a MILP, with an attack graph representing the cyber network and a state-space LTI model representing the physical system. Our novel objective function, based on MTTF, evaluates the impact of cyberattacks on the physical system. This framework serves as a tool for developing effective detection schemes and mitigation strategies in ICSs, especially in the absence of real-case attack data, by identifying various worst-case attack strategies. The experimental results validate our framework's efficiency and applicability. The identified worst-case cyberattacks can have more devastating physical consequences than a random cyberattack, thus directly highlighting its effectiveness in developing mitigating strategies. Nevertheless, our current proposed framework needs to include uncertainties present in both cyber and physical systems. The stochastic nature of these systems will be incorporated into our future analysis.

\normalsize
\bibliography{references}


\end{document}